\documentclass[psfig]{article}
\input{psfig}
\newtheorem{theo}{Theorem}
\newtheorem{lem}{Lemma}

\newtheorem{defi}{Definition}
\newtheorem{conj}{Conjecture}
\newtheorem{prop}{Proposition}

\newtheorem{q}{Question}
\newtheorem{cor}{Corollary}
\newtheorem{rem}{Remark}

\begin{document}

\makeatletter	   
\renewcommand{\ps@plain}{%
     \renewcommand{\@oddhead}{\textrm{Space of K\"{a}hler metrics}\hfil\textrm{\thepage}}%
     \renewcommand{\@evenhead}{\@oddhead}%
     \renewcommand{\@oddfoot}{}
     \renewcommand{\@evenfoot}{\@oddfoot}}
     \renewcommand{\thefootnote}{\fnsymbol{footnote}}
\makeatother     


\title{The Space of K\"{a}hler metrics }
\author{ Xiuxiong Chen\footnote{Research was supported partially by NSF postdoctoral fellowship.}  }
\date{Revised on July, 1999} 
\pagestyle{plain}
\bibliographystyle{plain}

\maketitle
\begin{abstract}  
Donaldson  conjectured\cite{Dona96}  that the space
of K\"ahler metrics  is geodesic convex by smooth geodesic
 and that it is a metric space. Following
Donaldson's program, we verify the second
part of Donaldson's conjecture completely and verify his first
part  partially. We also prove that the constant scalar curvature
metric is unique in {\bf each} K\"{a}hler class if the first Chern
class is either strictly negative or 0. Furthermore, if $C_1 \leq 0,$
the constant scalar curvature metric realizes the global minimum of
Mabuchi energy functional; thus it provides a new obstruction for
the existence of constant curvature metric: if the infimum of
Mabuchi energy (taken over all metrics in a fixed K\"{a}hler class)
isn't bounded from below, then there doesn't exist a constant curvature
metric.    This  extends the work of Mabuchi and Bando\cite{Bando87}: 
they showed that Mabuchi energy bounded from below is a necessary condition for the
existence of K\"{a}hler-Einstein metrics in the first
Chern class.
\end{abstract}
\section{Introduction to the problem}
\subsection{Brief introduction to the classical
problems in K\"{a}hler geometry}
Let $V$ be a K\"{a}hler manifold.
E. Calabi conjectured in 1954 
that any (1,1) form  which represents $C_1(V)$ (the first Chern class) 
is the Ricci form of some K\"{a}hler metrics on $V.\;$ Yau~\cite{Yau78}, 
in 1978,  proved this Calabi's conjecture. Around the same time,
Aubin~\cite{Au76} and Yau proved independently  the existence 
of a K\"{a}hler-Einstein metric on a K\"{a}hler manifold with
negative  first Chern class (also a conjecture of E. Calabi). 
G. Tian~\cite{Tian87}, in 1987,
 proved the existence of K\"{a}hler-Einstein
metric in a canonical K\"{a}hler class on complex surfaces  if 
the first Chern class is positive and the group of automorphism is reductive.
For further references on this subject, see \cite{Tian96} and \cite{Tian97}. 
An important conjecture by Yau \cite{Yau92}
relates the existence
of K\"{a}hler-Einstein metrics to the stability in the sense of Hilbert
schemes and Geometric invariant theory. \\

 K\"{a}hler-Einstein metrics could be treated as a
special kind of extremal K\"{a}hler metrics.  The question
of extremal k\"{a}hler metric 
 was first raised by E. Calabi
in his paper\cite{calabi82} : he considered $L^2$ norm of curvature as a functional
from a given K\"{a}hler class; a critical point of this functional
is called  an ``extremal K\"{a}hler metric.'' 
He  showed that any extremal K\"{a}hler metric 
must be symmetric under a maximal compact subgroup of the
holomorphic transformation group.  Using this structure
theorem of Calabi,
Marc Levine \cite{Levin85}  was able to construct a
K\"{a}hler surface on which there is no extremal K\"{a}hler metric.
In 1992, D. Burns and P. de Bartolomeis \cite{BurnsBa92} also produced an
example
of non-existence of extremal K\"{a}hler metric; their example
suggests some new obstruction for the existence of extremal metrics 
which is related to  some borderline 
semi-stability of hermitian vector  bundle. LeBrun\cite{LeBrun95}
also demonstrated that the existence of critical K\"{a}hler metrics might
be tied up
with the stability of corresponding vector bundles.
Donaldson\cite{Dona96} thought  that Yau's conjecture \cite{Yau92} should be extend
over to the general extremal K\"{a}hler metrics.
For further references in the subject of extremal metrics,
 please see \cite{LeBrun-Simanca94}, \cite{LeBrun95}
\cite{futaki88}  and references therein.\\

 Futaki\cite{futaki83} in 1983 introduced an analytic
invariant for any K\"{a}hler manifold with positive first Chern class. The vanishing of this invariant is a necessary
condition for the existence of a K\"{a}hler-Einstein metric on the
manifold. Later Futaki and Calabi~\cite{calabi85} generalized the invariant
to any compact K\"{a}hler class. This generalized Futaki invariant,
i.e., Calabi-Futaki invariant, is an analytic obstruction to
the existence of constant scalar curvature metric in a 
K\"{a}hler manifold.
In the same paper, Calabi also shows that  constant scalar curvature
metric and  extremal K\"{a}hler metric with non-constant scalar curvature
do not co-exist in a single K\"{a}hler
class.\\ 

 For the uniqueness, the  known results are
as follows:
  1)in 1950s, E. Calabi showed the uniqueness of 
K\"{a}hler-Einstein metrics if $C_1 \leq 0.\;$
  2)in 1987, Mabuchi and S. Bando~\cite{Bando87} showed the uniqueness
of K\"{a}hler - Einstein metrics up to holomorphic transformation
if the first Chern class is positive. Recently, Tian and X.H. Zhu 
\cite{TianZhu98}
proved the uniqueness of K\"{a}hler-Ricci Soliton with respect to
a fixed holomorphic vector field on any K\"{a}hler manifolds with positive
first Chern class.
 Although  very little was known about the uniqueness of general
extremal K\"{a}hler metrics, most experts 
in K\"{a}hler geometry  expect that  the extremal K\"{a}hler metric is unique in each K\"{a}hler 
class up to holomorphic transformation. In \cite{chen943} (also
see  \cite{chen981} for further references), we
demonstrated two degenerate extremal K\"{a}hler metrics in the same
K\"{a}hler class with different energy levels and different symmetry
 groups: one example is due to Calabi, the other is due to the author. 
To my knowledge, it appears
that this is the only non-uniqueness example known today.\\

\noindent {\bf Main results.}  Mabuchi (\cite{Ma87})\footnote{Around
the same time with Mabuchi's work, Bourguignon J. P. has worked
on something similar in a related subject \cite{Bourg85}.}
 in 1987 defined a Riemannian metric  on the space of K\"ahler metrics,
under which it  becomes (formally) a non-positive curved infinite dimensional
 symmetric space. Apparently unaware of Mabuchi's work, 
Semmes \cite{Semmes92}  and Donaldson
 \cite{Dona96}  re-discover this same metric again from different 
angles.
In \cite{Semmes92}, Semmes S. first pointed out that the geodesic
equation is a homogeneous complex Monge-Ampere equation on a manifold
of one dimension higher. In
 \cite{Dona96}, 
Donaldson  further conjectured
that the space is geodesically convex and it is a genuine metric
space. We prove that it is at least convex by $C^{1,1}$\footnote{Here we mean the mixed second derivatives is uniformly bounded. See theorem 3 in section 3 
for details.}  geodesics,
and from which we conclude that the space is indeed a metric
space, thus verifying the second part of Donaldson's conjecture. Moreover,
this $C^{1,1}$ geodesic realizes the absolute minimum of length
over all possible paths connecting the end points; thus the 
metric aforementioned is a genuine one.
Using these results, we are able to show that the constant
curvature metric is unique in each K\"{a}hler class if $C_1 < 0$
or $C_1 = 0.\;$ 
 Furthermore, if $C_1 \leq 0,$ we show that constant
scalar metric (if exists) realizes the global minimum of 
Mabuchi energy, which gives an affirmative answer to a question
raised by Gang Tian \cite{TianPrivate2} in this special case. This last statement also extends the work of Mabuchi and Bando\cite{Bando87}: they showed that
Mabuchi energy bounded from below is a necessary condition for the
existence of K\"{a}hler-Einstein metrics in the first
Chern class. In the light of Tian's work\cite{Tian97} in which he shows that
in K\"hler manifold with positive first Chern class and no non-trivial
holomorphic fields, the K\"hler-Einstein metric exists if and
only if the Mabuchi functional is proper (he actually uses
an equivalent functional instead of the Mabuchi functional)\footnote{The sufficient
part of this result was proved in \cite{DingTian93}.}. One would like to ask:
is this still true for constant scalar curvature metrics\footnote{Tian inform
us that he\cite{tian98} has conjectured that constant scalar curvature metrics
exist if and only if Mabuchi functional is proper.}?\\

\noindent{\bf Organization}: In section 2, we first summarize
the different approaches taken by Mabuchi, Semmes and Donaldson
independently in the
space of K\"{a}hler metrics; then
we introduce this
Riemannian metric on this infinite dimensional space and prove that 
it has non-positive sectional curvature in the formal sense. Then we introduce  Donaldson's
two conjectures and reduce the 1st conjecture to the existence
problem for the complex homogeneous  Monge-Ampere equation with Drichelet boundary
data. Readers are alerted that material in section 2.3-2.5 is  essentially
 a  re-presentation of Donaldson's work \cite{Dona96}, 
 included here for the convenience of readers. In section 3, we prove that this geodesic(CHMA) equation always has a 
$C^{1,1}$ solution. In section 4, 
we prove that  a continuous solution to the geodesic(CHMA) equation  in some appropriate weak sense is
 unique. In section 5, we show that the geodesic distance defined
by the length of $C^{1,1}$ geodesic  satisfies  
the triangular inequality. Using this,  we prove the
space of K\"{a}hler metrics is a metric space. In section 6
we show that extremal K\"{a}hler metric is unique in each K\"{a}hler
class if either $C_1(V) < 0$ or $C_1(V) = 0.\;$\\

\noindent{\bf Acknowledgment}: The author is very grateful for
Simon Donaldson who not only introduced him to this problem,
but also spent many long sessions with  him. He is also grateful
for the constant encouragement and support of E. Calabi, L. Simon
and R. Schoen during this work. he also want to thank Professor
 E. Stein for his help in soblev functions and embedding theorems.
 The author wish to thank S. Semmes
for points out some important references in this problem. Thanks also
goes to W.Y. Ding and his colleagues, and P. Guan for pointing out
some errors in an earlier version of this paper.

\section{Space of K\"{a}hler metrics}
\subsection  {Mabuchi and S. Semmes' Ideas}
Shortly after introducing the now famous Mabuchi functional, 
Mabuchi \cite{Ma87} set out and
defined a Riemannian metric in the space of K\"ahler metrics.
Besides showing formally it is a locally symmetric space with non-positive
sectional curvature, he also pointed out that the Mabuchi
energy is formally convex in this infinite dimensional space (in the
sense that the Hessian is semi-positive definite). Perhaps, this is
his original motivation for introducing such a metrics. Unaware
of Mabuchi's work,  in a remarkable paper \cite{Semmes92},
S. Semmes
studied the geometry of solution of complex Homogeneous Monge-Ampere
equation (CHMA). He observed that in some special domain $\Omega \times D$
where $\Omega$ is n-dimensional domain in $C^n$ and $D$ is a domain in complex
plane, the solution to CHMA is some sort of geodesic equation if the
data is rotationally symmetric when restricted to $D.\;$ He then
considered the space of pluri-subharmonic functions in $\Omega$
and defined a Riemannian metric in this space according to this
geodesic equation. It turns out that this space becomes non-positively
curved (locally) symmetric space in some formal sense.  
He also went on to study the variational problem of finding
a geodesic. It seems that he is mainly motivated from
providing a proper geometric meaning to solution of CHMA with right
 domain setting. 
Unlike real homogeneous Monge-Ampere equation (RHMA) whose solution always 
has proper geometric meaning, solution of a CHMA
 equation doesn't have a preferred geometric interpretation.
Without a proper geometry interpretation, it is very hard to work
on this subject. Of course, great progress has been made since
the famous work of L. Caffarelli, L. Nirenberg and J. Spruck\cite{CNS84}
and later their joint work with J. Kohn \cite{CKNS85}.  For instance, 
Lempert L. \cite{Lempt83}, E. Bedford and B.A. Taylor \cite{Bedford76};
leong P.\cite{Lelong86} and  important work of Krylov \cite{Krylov87}
and Evans \cite{Evans82}
$\cdots, $etc.. This is by no means a complete list of papers
in complex Monge-Ampere equation since the author is quite new to this
important field. For a complete and updated
references, please see S. Kolodziej \cite{Kolo98}.
Donaldson's recent work certainly makes Mabuchi and Semmes's original work
all the more remarkable.\\
 
\subsection { Brief summary of Donaldson's theory on space of K\"{a}hler
metrics}
 Motivated from complete different reasons, S. K. Donaldson  \cite{Dona96}
re-discovered this metric. More importantly, he outlined a strategy 
in \cite{Dona96} to relate this
geometry of infinite dimensional space to the existence problems in K\"ahler
geometry. In particular, he explains how one can uses this
extra structure in the infinite dimensional space to solve the problems of
the existence and uniqueness of extremal K\"ahler metrics.
In general, the later are  intractable problems from traditional means.
He regards the space of K\"{a}hler metrics in a fixed K\"{a}hler 
class as an infinite dimensional symplectic manifold with the automorphism group  $SDiff(V)$ (symplectic diffeomorphism group of $V$ into itself).
In \cite{Dona97}, he pointed out that  scalar curvature
 is the moment map
$\mu$ from this infinite dimension symplectic manifold to the
 dual space of the Lie algebra of its automorphism group\footnote{The {\it moment map} point of view here was also observed by A. Fujiki\cite{fujiki92}.}. Thus, to find an extremal K\"{a}hler metric
in a fixed K\"{a}hler class
in classical K\"{a}hler geometry  could be re-interpreted as
to  find a pre-image of $0$ of the moment map $\mu$ in this symplectic
setting. This acute observation sheds new
light into the otherwise intractable problem
of the existence of extremal K\"{a}hler metrics in a K\"{a}hler manifold;
at least conceptually, the picture looks much clear.
 He then proposed several conjectures whose ultimate resolution will
lead to a better understanding of extremal K\"{a}hler metric, and
for that matter, better understanding of K\"{a}hler geometry as well.
The most fundamental one among his conjectures is 
the so called geodesic conjecture:  any two
K\"{a}hler metrics in the same class is connected by a
smooth geodesic.
A second conjecture by him is that this
space  of K\"{a}hler metric is a metric space under this metric. If the geodesic
conjecture is true, this second conjecture will be a direct
consequence (since this space of K\"{a}hler metrics in a fixed
K\"{a}hler class  is non-positively curved in the formal sense.). He went 
on to show that the uniqueness
of extremal K\"{a}hler metric is a consequence of this
geodesic conjecture as well.\\

\subsection{ Riemannian metrics in the infinite dimensional space.}
 Now  we introduce this metric here. Readers are 
referred to Mabuchi,  S. Semmes and  Donaldson's original writing 
for details.
Consider the space of K\"{a}hler potentials in a fixed K\"{a}hler
class as:
\[
 {\cal{H}} = \{ \varphi \in C^{\infty}(V) : \omega_{\varphi} =  \omega_{0} + \sqrt{-1} \partial \overline{\partial} \varphi > 0 \; {\rm on}\; V\}.
\]
Clearly, the tangent space  $T \cal {H} $  is $C^{\infty} (V).\;$
Each K\"{a}hler potential  $\phi \in \cal {H}$ defines a measure
$d\,\mu_{\phi} = {1\over {n!}} \omega_{\phi}^n.\; $
Now we define a Riemannian metric on the infinite dimensional 
manifold $\cal {H}$ using the $L^2 $ norm provided by these measures.
A tangent vector in $\cal {H}$ is just a function
in $V.\;$ For any  vector
$\psi \in T_{\varphi} \cal {H}, $ we define the length of this
vector as
\[
\|\psi\|^2_{\varphi} =\int_{V}\psi^2\;d\;\mu_{\varphi}. 
\]  

For a path $\varphi(t) \in {\cal {H}} (0\leq t \leq 1),$
the length  is given by
\[
  \int_0^1 \sqrt{\int_V {\varphi(t)'}^2 d\,\mu_{\varphi(t)}} \; d\,t
\]

and the geodesic equation is
\begin{equation}
  \varphi(t)'' - {1\over 2} |\nabla \varphi'(t)|^2_{ \varphi(t)} = 0,
\label{geodesic}
\end{equation}
where the derivative and norm in the 2nd term of the left hand side
are taken with respect to the metric $\omega_{\varphi(t)}.\;$\\

This geodesic equation shows us how to define a connection on the tangent 
bundle of ${\cal H}$. The notation is simplest if one thinks of such a connection 
as a way of differentiating vector fields along paths. Thus, if $\phi(t)$ is 
any path in ${\cal H}$ and $\psi(t)$ is a field of tangent vectors 
along the path 
(that is, a function on $V \times [0,1]$), we define the covariant derivative 
along the path to be
$$ D_{t}\psi = \frac{\partial\psi}{\partial t} - {1\over 2}  (\nabla \psi, \nabla 
\phi')_{\phi}.  $$
This connection is torsion-free because in the canonical \lq\lq co-ordinate 
chart'', which represents ${\cal H}$ as an open subset of $C^{\infty}(V)$, the 
\lq\lq Christoffel symbol''
$$\Gamma: C^{\infty}(V) \times C^{\infty}(V) \rightarrow C^{\infty}(V)$$ at 
$\phi$
is just
$$\Gamma(\psi_{1}, \psi_{2}) = - {1\over 2}  (\nabla \psi_{1}, 
\nabla\psi_{2})_{\phi}$$
which is symmetric in $\psi_{1},\psi_{2}$. The connection is metric-compatible 
because
$$\begin{array}{lrl} {1\over 2} \frac{d}{dt} \Vert \psi\Vert^{2}_{\phi} & = &\frac{d}{dt}\int_{V}\psi^{2} 
d\mu_{\phi}\\
&= & \int_{V}  \frac{\partial \psi}{\partial t} \psi + {1\over 2}  \psi^{2} \Delta 
(\phi') d\mu_{\phi}\\
&= & \int_{V} \frac{\partial \psi}{\partial t} \psi - {1\over 4}(\nabla(\psi^{2}), 
\nabla \phi')_{\phi}\  d\mu_{\phi}\\
&= &  \int_{V} (\frac{\partial \psi}{\partial t}  - {1\over 2} (\nabla\psi, \nabla 
\phi')_{\phi}\ ) \psi d\mu_{\phi}\\
&= & \langle D_{t} \psi, \psi\rangle.\end{array}$$
Here $\triangle$ is complex Laplacian operator. The main theorem 
proved  in \cite{Ma87}(and later reproved in \cite{Semmes92} and
 \cite{Dona96}) is:\\
 
\noindent {\bf Theorem A} {\it The Riemannian manifold $\cal {H} $ is an infinite dimensional 
symmetric space; it admits a Levi-Civita connection whose curvature is 
covariant constant. At 
 a point 
$\phi\in{\cal {H}}$ the curvature  is given by
$$      R_{\phi}(\delta_{1}\phi, \delta_{2}\phi) \delta_{3}\phi=
- {1\over 4}  \{ \{ \delta_{1}\phi, \delta_{2}\phi\}_{\phi}, 
\delta_{3}\phi\}_{\phi},$$
where $\{\ ,\ \}_{\phi}$ is the Poisson bracket on $C^{\infty}(V)$ of the 
symplectic form $\omega_{\phi}$; and $\delta_1 \phi, \delta_2 \phi \in T_{\phi} {\cal H}.$ }

(Recall that in infinite dimensions the usual argument gives the uniqueness of 
a Levi-Civita [i.e. torsion-free, metric-compatible] connection, but not the 
existence in general.)
        The formula for the curvature of ${\cal {H}}$ entails that the sectional 
curvature is non-positive, given by
$$    K_{\phi}(\delta_{1}\phi, \delta_{2}\phi) = - {1\over 4}  \Vert \{ \delta_{1} 
\phi , \delta_{2}\phi\}_{\phi}\Vert_{\phi}^2. $$
 
Different proofs of this theorem have been appeared in  \cite{Ma87}, \cite{Semmes92} and \cite{Dona96}. We will skip the proof here, interested readers are 
referred to these papers if they are interested in the proof.\\

The expression for the curvature tensor in terms of Poisson brackets shows 
that $R$ is invariant under the action of the symplectic-morphism group. Since the 
connection on $T{\cal H}$ is induced from an ${\rm SDiff}$-connection, it follows 
that $R$ is covariant constant, and hence ${\cal H}$ is indeed an 
infinite-dimensional symmetric space.

\subsection{Splitting of ${\cal H}$}
 There is obviously a decomposition of the tangent space:
$$     T_{\phi}{\cal H} = \{ \psi : \int_{V} \psi d\mu_{\phi}=0 \} \oplus {\bf R}. 
$$
We claim that this corresponds to a Riemannian decomposition 
$${\cal H} = {\cal H}_{0}\times {\bf R}. $$
We are interested to see this Riemannian splitting  more explicitly, 
partly because we see the appearance 
of a 
functional $I$ on the space of K\"{a}hler potentials, which is well-known
 in the literature, see \cite{Au84}, \cite{Tian97} for example. 
The decomposition of tangent space of $\cal H$ gives a $1$-form
$\alpha$ on ${\cal H}$ with
$$  \alpha_{\phi}(\psi) = \int_{V} \psi d\mu_{\phi}, $$
and it is straightforward to verify that this $1$-form is {\it closed}. 
Indeed

$$ (d\alpha)_{\phi}(\psi,\tilde{\psi}) =  \int_{V} \left(\tilde{\psi}\Delta\psi- \psi 
\Delta\tilde{\psi}\right)  = 0.\;$$
This means that there is a function
$I:{\cal H}\rightarrow {\bf R}$ with $I(0)=0$ and $dI=\alpha$, and it is this function 
which gives rise to the corresponding Riemannian decomposition. We call a K\"{a}hler 
potential $\phi$ {\it normalized} if $I(\phi)=0$. Then any K\"{a}hler metric has a 
unique normalized potential, and the restriction of our metric on ${\cal H}$ to 
$I^{-1}(0)$ endows the space  ${\cal H}_{0}$ of K\"{a}hler metrics with a Riemannian 
structure; this is 
independent of the choice of base point $\omega_{0}$ and clearly makes 
${\cal H}_{0}$ into a  symmetric space. 
The functional $I$ can be written more explicitly by integrating $\alpha$ 
along lines in ${\cal H} $ to give the formula
$$   I(\phi) =\sum_{p=0}^{n} \frac{1}{(p+1)! (n-p)!} \int_{V} \omega_{0}^{n-p} 
(\partial \overline{\partial}\phi)^{p}\ \phi           
.$$

\subsection{Donaldson' Conjectures}
We will now study the geodesic equation  in $\cal {H}$ in more detail, and 
interpret the solutions geometrically. Suppose $\phi_{t},\ t\in [0,1]$, is a
path in $\cal {H}$. We can view this as a function on $V\times [0,1]$ and in turn 
as a function on $V\times [0,1]\times S^{1}$, with trivial dependence on the 
$S^{1}$ factor; that is, we define
                $$   \Phi(v,t,e^{is}) = \phi_{t}(v). $$
We regard the cylinder ${\bf R}=[0,1]\times S^{1}$ as a Riemann surface with 
boundary in the standard way---so $t+is$ is a local complex co-ordinate. Let
$\Omega_{0}$ be the pull-back of $\omega_{0}$ to $V\times {\bf R}$ under the 
projection map and put $\Omega_{\Phi} = \Omega_{0}+ \partial \overline{\partial} \Phi$, a $(1,1)$-form 
on $V\times {\bf R}$. 
Then we have:\\

\begin{prop} The path $\phi_{t}$ satisfies the geodesic equation (\ref{geodesic}) if and only if 
$\Omega_{\Phi}^{n+1} =0$ on $V\times {\bf R}$.
\end{prop}

\noindent {\bf Proof:} Denote the metric defined by $\omega_0, \omega_{\phi}$ as $g, g'.\;$ Then
$$ {1\over {n!}} \omega_{\phi}^n = \det\, g';\qquad {1\over {n!}} \omega_{0}^n = \det\, g.
$$ 
Then geodesic equation is equivalent to the following
(if $det \;g' \neq 0$)
\[ 
  (\phi'' - {1\over 2} \mid \nabla \phi'\mid_{g'}^2 ) \;det\; g' = 0.
\]
The last equation is equivalent to
\[
det  \left( \begin{array}{cc} g' & \left( \begin{array} {c} {{\partial \phi'} \over {\partial z_1 }} \\
{{\partial \phi'} \over {\partial z_2 }}\\
\downarrow\\{{\partial \phi'} \over {\partial z_n }} \end{array}\right)\\
         \left( \begin{array} {cccc} {{\partial \phi'} \over {\partial \overline{z_1} }} &
{{\partial \phi'} \over {\partial \overline{z_2} }} &
\cdots & {{\partial \phi'} \over {\partial \overline {z_n} }} \end{array}\right)              &  \phi''   \end{array} \right) = 0.
\]
Let $w = t + \sqrt{-1} s,$ then $t = Re(w).\; $ The above equation
could be re-written as
\[
det  \left( \begin{array}{cc} (g +{{\partial^2 \phi}\over{\partial {z_{\alpha}} \partial {\overline z_{\beta}}}}  )_{n\,n} & \left( \begin{array} {c} {{\partial^2 \phi} \over {\partial z_1 \partial {\overline w} }} \\
{{\partial^2 \phi} \over {\partial z_2 \partial {\overline w} }}\\
\downarrow\\{{\partial^2 \phi} \over {\partial z_n \partial {\overline w}}} \end{array}\right)\\
         \left( \begin{array} {cccc} {{\partial^2 \phi} \over {\partial \overline{z_1} \partial {w}  }} &
{{\partial^2 \phi} \over {\partial \overline{z_2} \partial {w} }} &
\cdots & {{\partial^2 \phi} \over {\partial \overline {z_n} \partial {w}}} \end{array}\right)              &  {{\partial^2 \phi}\over{\partial {w} \partial {\overline w}}}   \end{array} \right) = 0.
\] 
This is just $\Omega_{\Phi}^{n+1} =0.\;$ The proposition is then
proved.  $\qquad QED.$ \\

 Given boundary data ---a real value function $\rho\in C^{\infty} (\partial 
(V \times {\bf R})),$ we consider the set of functions $\Phi$ on $V\times {\bf R}$ 
which agree with $\rho$ on the boundary. Then we define the 
variation of $I_{\rho}$ on this set by
$$ \delta I_{\rho} = \frac{1}{(n+1)!}\int_{V\times {\bf R}} \delta \Phi \ 
\Omega_{\Phi}^{n+1} ,  $$
where the variation $\delta \Phi$ vanishes on the boundary by hypothesis. This 
boundary condition means that we can  show easily 
that this formula defines a functional $I_{\rho}.\;$ To prove this,
one only need to show that the second derivatives of $I_{\rho}$ with
respect to two infinitesimal variation $\delta_1 \Phi$ and $ \delta_2 \Phi $
is symmetric. The second derivatives is:
\[ {1\over 2} \cdot \frac{1}{(n+1)!}
\int_{V} \delta_1 \Phi \; \triangle \;\delta_2 \Phi \;\Omega_{\Phi}^{n+1}
\] 
which is clearly symmetric on $\delta_1 \Phi$ and $ \delta_2 \Phi.\; $
Here $\triangle$ is the Laplacian operator of $\Omega_{\Phi}$
on $V \times {\bf R}.\;$\\

 This functional $I_{\rho}$ 
reduces to the energy functional on paths, by an integration 
by parts, in the case when
 ${\bf R}$ is the cylinder  and we restrict to $S^{1}$-invariant data.
Suppose $\phi(t) (0\leq t \leq  1)$ is a path in $\cal H, \;$
and $\delta \phi$ represents the infinestimal variation of $\phi$
while keep value of $\phi$ fixed when $t=0, 1.\;$ Thus, the variation
of $I_{\rho}$ in $\delta \phi$ direction is (follow notations in the
proof of previous proposition):
\[
  \delta I_{\rho} = \frac{1}{(n+1)!}\int_{V\times R} \delta \phi\; 
\Omega_{\Phi}^{n+1} = \frac{1}{(n+1)!}\int_{t=0}^{1} \,\int_{V} \delta \phi 
(\phi''-   {1\over 2} \mid \nabla \phi'\mid_{g'}^2 ) \;det\; g' d\,t.  
\]
 On the other hand, the variation of energy functional along this path
is:
\[
\delta E = \int_{t=0}^{1}\; \int_V \delta \phi\;
(\phi''-   {1\over 2} \mid \nabla \phi'\mid_{g'}^2 ) \;det\; g' d\,t
\]
where $E =  \int_{t=0}^{1}\; \int_V \phi'(t)^2 \;det\; g' d\,t.\;$
Thus, in  case when
 ${\bf R}$ is the cylinder  and we restrict to $S^{1}$-invariant data,
$I_{\rho}$ equal to the energy functional on the path up to a multiple
of constant.\\

The following is the first conjecture by Donaldson in ~\cite{Dona96}:

\begin{conj}(Donaldson)
 Let ${\bf R}$ be a compact Riemann surface with boundary and $\rho:V\times 
\partial R\rightarrow {\bf R}$ be a function such that $\omega_{0}- {\sqrt{-1}\;\overline{\partial}\partial} \rho$ is 
a strictly positive (1,1) form on each slice $V\times \{z\}$ for each fixed 
$z\in \partial R$. Let $\cal{S}_{\rho}$ be the set of functions $\Phi$ on $V\times 
R$ equal to $\rho$ over the boundary and such that
 $\omega_{0} -{\sqrt{-1}\;\overline{\partial}\partial}\Phi$ is  strictly positive on  every 
slice $V\times \{w\}, w\in R$. Then there is a unique solution of the 
Monge-Ampere equation $(\Omega_{0} -{\sqrt{-1}\;\overline{\partial}\partial}\Phi)^{n+1}=0$ 
in $\cal{S}_{\rho}$, and this solution realizes the absolute minimum of 
the functional $I_{\rho}$.
\end{conj}

This  question  is a version of the Dirichlet problem for the 
complete degenerate Monge-Ampere equation, a topic around which there is a substantial
literature; see \cite{Au84},\cite{Klimek91} for example. Note that regularity questions are 
very important in this theory, since the equation is not elliptic.\\

 In the case of the geodesic problem, when 
the functional can be rewritten as the  energy of a path; if these 
infimum are strictly positive, 
for all choices of fixed, distinct,  end points, they make ${\cal H}$ into a metric space, in the usual fashion. In this 
connection,  Donaldson proposes the following conjecture (after
verifying that it will be satisfied by a smooth geodesic):

\begin{conj}(Donaldson)
If $\phi\in{\cal {H}}_{0}$ is normalized and $\tilde{\phi}_{t},\ t\in[0,1]$
 is { any} path from $0$ to $\phi$ in ${\cal {H}}$ then
\begin{equation}  \int_{0}^{1} \int_{V} \left(\frac{d\tilde{\phi}}{dt}\right)^{2}  
d\mu_{\tilde{\phi}_{t}} dt \geq M^{-1} \left( \max( \int_{\phi>0} \phi 
d\mu_{\phi}, -\int_{\phi<0} \phi d\mu_{0}) \right)^{2}. \label{eq:lowerbound}
\end{equation}
\end{conj}

The restriction to normalized potentials  $\phi$ is not important since we know that ${\cal {H}}$ 
splits as a product, and we could immediately write down a corresponding 
inequality, involving $I(\phi)$, for any $\phi\in{\cal {H}}$. If this
conjecture and the geodesic conjecture are proved, then $\cal H$ is a metric space.

we want to use continuous method to treat this existence problem of geodesic 
between any two points in $\cal H$.

\section {Existence of $ C^{1,1} $ solution}
Let $V$ be a $n-$ dimensional K\"{a}hler manifold without boundary, ${\bf R}$ be a
Riemann surface with boundary. The case we concerned most is
when ${\bf R}$ is a cylinder. Suppose $g = g_{\alpha \overline{\beta}} dz_{\alpha} d\,\overline{z_{\beta}} (1\leq \alpha,\beta \leq n)$ is a given K\"{a}hler metric in $V.\;$ Then $ \tilde{g} =  g_{\alpha \overline{\beta}} dz_{\alpha} d\,\overline{z}_{\beta} + dw\, \overline{dw}$ is
a K\"{a}hler metric in $V\times {\bf R},\;$ and
$\tilde{\varphi} = \varphi -  |w|^2.\;$ 
For convenience, we still denote $\tilde{g}$ as $g$, and $\tilde{\varphi}$
as $\varphi$  when there is no confusion arisen. 
Also, let $z_{n+1}=w.\;$ Then
$z= (z_1,z_2,\cdots,z_n,z_{n+1})$ is a point in
$V\times {\bf R} $ and $z'=(z_1, z_2,\cdots z_n)$ is a
point in $V.\;$ Let $\varphi(z) = \varphi(z', w)$ be a function
in $V\times {\bf R}$ such that $g + \partial_{z'} \overline {\partial_{z'}} 
\varphi(z',w)$ is a K\"{a}hler metric in $V$ for each $w \in {\bf R}.\;$
We want to solve the degenerated Monge-Ampere equation:
\begin{equation}
   det\;(g + {{\partial^2 \varphi}\over{\partial z_{\alpha}  \partial \overline{z}_{\beta}}})_{(n+1)(n+1)}  =0 \;\; {\rm in}\; V\times {\bf R}; \qquad {\rm and}\;\;
\varphi = \varphi_0 \;{\rm in}\;\; \partial (V\times {\bf R}). 
\label{eq:euler1}
\end{equation}

  We want to use the continue method to solve this equation. 
Consider the continuous equation $0 \leq t \leq 1.\;$
\begin{equation}
   det\;(g + {{\partial^2 \varphi}\over{\partial z_{\alpha}  \partial \overline{z}_{\beta}}}) = t \;  det\;(g + {{\partial^2 \varphi_0}\over{\partial z_{\alpha}  \partial \overline{z}_{\beta}}}), \;\; {\rm in}\; V\times {\bf R}; \qquad {\rm and}\;\;
\varphi = \varphi_0 \;{\rm in}\;\; \partial (V\times {\bf R}). 
\label{eq:euler2}
\end{equation}
Suppose $\varphi_0$ is a solution to (\ref{eq:euler2}) at $t=1$  such that $\displaystyle \sum_{\alpha,\beta = 1}^{n+1} (g_{\alpha \overline{\beta}} + {{\partial^2 \varphi_0}\over{\partial z_{\alpha}  \partial \overline{z}_{\beta}}}) dz_{\alpha} d\,\overline{z}_{\beta} $
is strictly positive K\"{a}hler metric in $V \times {\bf R}$\footnote{By definition, for any $\varphi_0 \in {\cal H},$ 
 $\displaystyle \sum_{\alpha,\beta = 1}^{n+1} (g_{\alpha \overline{\beta}} + {{\partial^2 \varphi_0}\over{\partial z_{\alpha}  \partial \overline{z}_{\beta}}}) dz_{\alpha} d\,\overline{z}_{\beta} $
is strictly positive K\"{a}hler metric in each $V-$ slice $V \times \{w\}.\;$ Let $\Psi$ be a strictly convex function of $w$ which vanishes on $\partial {\bf R}.\;$ Then for large enough constants $m , \displaystyle \sum_{\alpha,\beta = 1}^{n+1} (g_{\alpha \overline{\beta}} + {{\partial^2 (\varphi_0 + m \Psi) }\over{\partial z_{\alpha}  \partial \overline{z}_{\beta}}}) dz_{\alpha} d\,\overline{z}_{\beta} $
is a strictly positive K\"{a}hler metric in $V \times {\bf R}.\;$   }. Denote
$f =  det\;(g + {{\partial^2 \varphi_0}\over{\partial z_{\alpha}  \partial \overline{z}_{\beta}}}) (det \;g)^{-1} > 0.\;$ Then equation (\ref{eq:euler2}) can be re-written
in a better form
\begin{equation}
   det\;(g + {{\partial^2 \varphi}\over{\partial z_{\alpha}  \partial \overline{z}_{\beta}}}) = t \cdot f\cdot  det\;(g) \;\; {\rm in}\; V\times {\bf R}; \qquad {\rm and}\;\;
\varphi = \varphi_0 \;{\rm in}\;\; \partial (V\times {\bf R}). 
\label{eq:euler3}
\end{equation}
Clearly, $\varphi_0$ is the unique solution to this equation at $t=1.\;$
Since the equation is elliptic, this equation can be uniquely solved
for $t $ sufficiently closed to $1 $(the kernal of linearized operator is zero
for any $t > 0$). Let $t_0$ be such that
(\ref{eq:euler3}) has a unique smooth solution for every $t \in (t_0,1].\;$ We want
to show that $t_0=0 $ in this section. Observe that equation (\ref{eq:euler3}) is elliptic
for every $t > 0.\; $ Hence, the solution will be as smooth
as the boundary value once we show that 2nd derivatives of $\varphi$
is uniformly bounded. Let $h$ be a super harmonic function on $V\times {\bf R}$ 
with respect to $g$ such that $\triangle_g h + n + 1 = 0.$ and 
 $h = \varphi_0$ in $ \partial (V\times {\bf R} ).\;$
Then for any solution of equation (\ref{eq:euler3}) for $t < 1$, we have $C^0$
 bound of the solution:

\begin{lem} If $\varphi$ is a solution of equation (\ref{eq:euler3}) at $0<t <1,$
then $\varphi$ has the following a priori $C^0 $ estimate due to 
maximum principal:
\[   \varphi_0 \leq \varphi \leq h,\qquad {\rm in } \; V\times {\bf R}.
\] 
\end{lem}
 
 For $C^2$ estimate, we follow Yau's famous work in Calabi's conjecture.
Essentially, we reduce it to a boundary estimate since we have $C^0$ estimate:
\begin{lem}(Yau) If $\varphi$ is a solution of equation (\ref{eq:euler3}) at $0<t <1,$
then $\varphi$ has the following a priori $C^2 $ estimate:
\[
\begin{array}{lcl}
  \triangle' ( e^{- C\varphi} (n+1 + \triangle \varphi)) & \geq & e^{-C\varphi} (\triangle \ln f - (n+1)^2 \displaystyle \inf_{i\neq l} (R_{i\overline{i} l\overline{l}})) - C e^{-C\varphi} (n+1) (n+1 + \triangle \varphi) \\ & & 

+ (C +\displaystyle \inf_{i\neq l} (R_{i\overline{i} l\overline{l}})) e^{-C\varphi} (n+1 + \triangle \varphi)^{ 1 + {1\over n}} ( t f)^{-1}.
\end{array}
\]   
where $C + \displaystyle \inf_{i\neq l} (R_{i\overline{i} l\overline{l}}) > 1,\;$ $\triangle $ is the Laplacian operator with respect to to $g,$ while
 $\triangle' $ is the Laplacian operator with respect to to $g' = g + {{\partial^2 \varphi}\over{\partial z_{\alpha}  \partial \overline{z}_{\beta}}} d\,z_{\alpha}\; \overline{d\, z_{\beta}} $ and $R_{i\overline{i} l\overline{l}}$ is
the Riemannian curvature of $g.\;$
\end{lem}
From the a priori estimate in Lemma 2,  either
$e^{- C\varphi} (n+1 + \triangle \varphi)$ is uniformly bounded 
in $V \times {\bf R}$ or it achieves maximum value at $\partial (V \times {\bf R}).\;$
Lemma 1 asserts that $\varphi$ is uniformly bounded
from above and below, then 
\begin{cor} There exists a constant $C$ which depends only on
$(V\times {\bf R}, g) $ such that
\[
\displaystyle \max_{V \times {\bf R}}\; (n+1 + \triangle \varphi) 
\leq C ( 1 + \displaystyle \max_{\partial (V \times {\bf R})}\; (n+1 + \triangle \varphi)). \]
\end{cor}
\begin{theo} If $\varphi$ is a solution of equation (\ref{eq:euler3}) at $0<t <1,$
then  there exists a constant $C$ which depends only on
$(V\times {\bf R},g)$ such that:
\begin{equation}
\displaystyle \max_{V \times {\bf R}}\; (n+1 + \triangle \varphi) 
\leq C \displaystyle \max_{V \times {\bf R}} \;(|\nabla \varphi|_g^2+1).
\end{equation}
\end{theo}
In light of Corollary 1, we only need to prove the inequality (6)
on the boundary, i.e,,
\[
\displaystyle \max_{\partial (V \times {\bf R})}\; (n+1 + \triangle \varphi) 
\leq C \;  \displaystyle \max_{V \times {\bf R}} \;(|\nabla \varphi|_g^2 + 1).
\]
We will prove this inequality in the next subsection.

\begin{theo}  If $\varphi_i (i = 1, 2, \cdots)$ are  solutions of equation (\ref{eq:euler3}) at $0<t_i <1,$ and the inequality (6) holds uniformly for all these
solutions $\{\varphi_i, i \in {\bf N}\}$, then there exists a constant
 $C_1$ independent of $i$ such that
\[
   \displaystyle \max_{V \times {\bf R}}\; (n+1 + \triangle \varphi) 
\leq C \;  \displaystyle \max_{V \times {\bf R}} \;(|\nabla \varphi|_g^2 + 1)
  < C_1.
\]
\end{theo}
This is proved via a blowing up argument. We will show this in
subsection 3.2. 
\begin{rem} By now it is standard estimate of Monge-Ampere equations,
that if
 \[ \displaystyle \max_{V \times {\bf R}}\; (n+1 + \triangle \varphi) 
\leq C \;  \displaystyle \max_{V \times {\bf R})} \;(|\nabla \varphi|_g^2 + 1)
  < C_1 \]
then equation (\ref{eq:euler3}) for $ t_1, t_2 ,\cdots$ is a sequence
of uniform elliptic equations. The higher derivative of
the solution $\varphi_i$ has a uniform bound as long
as $\displaystyle \liminf_{i \rightarrow \infty} t_i > 0.\;$
\end{rem}
\begin{theo} There exists a $C^{1,1} (V\times {\bf R})
$ function which solves equation (\ref{eq:euler1})
weakly.  In other words,  for any two points
$\varphi_0, \varphi_1 \in \cal H,$ there exists a  geodesic
path $\varphi(t): [0,1] \rightarrow \overline{\cal H} $ and
a uniform constant $C$ such
that the following holds:
\[
0 \leq \left( g_{i \overline{j}}\; + \; {{\partial^2 \varphi} \over 
{\partial z_i \partial \overline{z_j}}} \right)_{(n+1) (n+1)} \leq C \left({\tilde{g}}_{i \overline{j}}\right)_{(n+1)(n+1)}.
 \]
Here $z_1, z_2 \cdots, z_n$ are
local coordinates in $V$ and $ t = Re \;(z_{n+1}).\;$ And $\tilde{g} =  g_{\alpha \overline{\beta}} dz_{\alpha} d\,\overline{z}_{\beta} + dw\, \overline{dw} $ is a
fixed product metric in $V \times {\bf R}.\;$

\end{theo}

Following notations in theorem 2, we want to show that $t_0 =
\displaystyle \liminf_{i \rightarrow \infty} t_i  =0.\; $ Otherwise,
assume $t_0>0.\;$ Then equation (\ref{eq:euler3}) has a unique smooth solution
for $1 \geq t > t_0.\;$ Following from theorem 2, then we have
uniform upper bound for $\triangle \varphi + (n+1)$ for all $t_i > t_0 > 0.\;$
Then equation (\ref{eq:euler3}) implies that $g'_i = g + {{\partial^2 \varphi_i}\over{\partial z_{\alpha}  \partial \overline{z}_{\beta}}} d\,z_{\alpha}\; \overline{d\, z_{\beta}} $ is bounded uniformly from 
below by a uniform positive constant (this positive low
bound approaches 0 when $t \rightarrow 0)$. Thus, from equation (\ref{eq:euler3}),
we obtain uniform higher derivative estimates for solution
$\varphi_i.\;$ Therefore these solution converge to a regular
solution at $t_0 > 0.\;$ Again, since equation (\ref{eq:euler3}) at $t_0$ is
an elliptic equation and the kernal of the linearized operator
is zero, it can then be solved for any $t$ sufficiently
closed to $t_0.\;$ But this contradicts to the definition of
$t_0.\;$ Thus $t_0 = 0.\;$ We can choose a subsequence
of $t_i \rightarrow 0$ such that $\varphi_i$ converge weakly
in $ C^{1,1} (V \times {\bf R}) $
where $\Omega$ is relative compact subset of $ V \times {\bf R}.$
Again via maximum principal, we can show this limit is
unique and define a weak solution of equation (\ref{eq:euler1}).
\subsection {Boundary estimate}
We want to estimate  $ \triangle \varphi $  at any point  in 
the boundary $ \partial (V \times {\bf R}) = V \times \partial {\bf R}.\;$ 
Let $p$
be a generic point in $\partial (V \times {\bf R}).\;$ Now choose
a small neighborhood $U$ of $p$ in $V\times {\bf R}$ (this
will be a half geodesic ball since $ p \in \partial ( V\times {\bf R}))$
and a local coordinate chart such that $ g_{\alpha \overline{\beta}}(p) = \delta_{\alpha \overline{\beta}}$ and $p = (z=0)$
\[
        {1\over 2}   \delta_{\alpha \overline{\beta}}   \leq g_{\alpha \overline{\beta}}(q)  \leq 2 \delta_{\alpha \overline{\beta}} ,\qquad \forall q\; \in\; U.\;
\]
Since $ \displaystyle \sum_{\alpha,\beta = 1}^{n+1}(g_{\alpha \overline{\beta}}  + {{\partial^2 \varphi_0}\over{\partial z_{\alpha}  \partial \overline{z}_{\beta}}}) d\,z_{\alpha} d\,\overline{z}_{\beta} $ is a positive K\"{a}hler metric
in $V\times {\bf R},$ there exists a constant $\epsilon > 0$ such
that
\[  g_{\alpha \overline{\beta}} + {{\partial^2 \varphi_0}\over{\partial z_{\alpha}  \partial \overline{z}_{\beta}}} > 2\;\epsilon \cdot g_{\alpha \overline{\beta}}, \qquad {\rm in}\; \; V \times {\bf R}.
\]
In the neighborhood $U$ of $p,$ we have
\begin{equation}  g_{\alpha \overline{\beta}} + {{\partial^2 \varphi_0}\over{\partial z_{\alpha}  \partial \overline{z}_{\beta}}} > \epsilon \cdot \delta_{\alpha \overline{\beta}}\qquad {\rm in}\; \; V \times {\bf R}.
\end{equation}
  We have the 
trivial estimates in $\partial (V \times {\bf R})$:
\[
{{ \partial (\varphi - \varphi_0)}\over {\partial z_{\alpha}}} = 0,\qquad
{{\partial^2 (\varphi-\varphi_0)}\over
{\partial z_{\alpha} \partial \overline{z}_{\beta}}} = 0,\qquad \forall
\; 1\leq \alpha, \beta \leq n.
\]
In order to estimate  $\triangle \varphi =  \displaystyle \sum_{\alpha,\beta=1 }^{n+1}g^{\alpha \overline{\beta}} \; {{\partial^2 \varphi}\over
{\partial z_{\alpha} \partial \overline{z}_{\beta}}} $ in $\partial (V \times {\bf R}),$ 
we only need to estimate ${{\partial^2 (\varphi-\varphi_0)}\over
{\partial z_{\alpha} \partial \overline{z}_{\beta}}} $ when either
$\alpha $ or $\beta$ is $n+1.\;$ We will estimate ${{\partial^2 (\varphi-\varphi_0)}\over
{\partial z_{\alpha} \partial \overline{z}_{n+1}}}  (\alpha \leq n)$ first,
then use equation (\ref{eq:euler3}) to derive estimate for ${{\partial^2 (\varphi-\varphi_0)}\over
{\partial z_{n+1} \partial \overline{z}_{n+1}}}.\; $\\

Now we set up some conventions: 
\[
  z_{\alpha} = x_{\alpha} + \sqrt{-1}\; y_{\alpha},\;\; \forall \; 1\leq \alpha \leq n;\qquad z_{n+1} = x + \sqrt{-1} \;y
\]
where $ {\bf R} $ near $\partial {\bf R}$ is given by $ x\geq 0$. 

\begin{lem} There exists a constant $C$ which depends only on $(V\times {\bf R}, g) $ such that
\[ 
 | {{\partial^2 \varphi }\over
{\partial z_{\alpha} \partial \overline{z}_{n+1}}}(p)| \leq C ( \displaystyle
\max_{V\times {\bf R}} \; |\nabla \varphi|_g + 1) .  
\] 
\end{lem}
{\bf Proof of theorem 1}: At point $p$, equation (\ref{eq:euler3}) reduces to
\[
  det ( \delta_{\alpha \overline{\beta}} + {{\partial^2 \varphi}\over
{\partial z_{\alpha} \partial \overline{z}_{\beta}}}) = t \cdot f.
\]
In other words,
\[
    {{\partial^2 \varphi}\over
{\partial z_{n+1} \partial \overline{z}_{n+1}}} = t \cdot f - {{\partial^2 \varphi}\over
{\partial z_{\alpha} \partial \overline{z}_{n+1}}}\cdot {{\partial^2 \varphi}\over
{\partial \overline{z}_{\alpha} \partial {z}_{n+1}}}.
\]
Lemma 3 then implies that
\[
    | {{\partial^2 \varphi}\over
{\partial z_{n+1} \partial \overline{z}_{n+1}}}| \leq C (\displaystyle
\max_{V\times {\bf R}}|\nabla \varphi|_g^2 + 1  ).
\]
Then,
\[
   |\triangle \varphi(p)| =   |\displaystyle \sum_{\alpha,\beta=1 }^{n+1}g^{\alpha \overline{\beta}} \; {{\partial^2 \varphi}\over
{\partial z_{\alpha} \partial \overline{z}_{\beta}}}(p)| \leq  C (\displaystyle
\max_{V\times {\bf R}}|\nabla \varphi|_g^2 + 1  ).
\]
Since $p$ is a generic point in $\partial (V \times {\bf R}),$ then
theorem 2 holds true.       QED.\\

Let $D$ be any constant linear 1st order operator near the boundary
( for instance $ D = \pm {\partial \over {\partial x_{\alpha}}},\; \pm 
{\partial \over {\partial y_{\alpha}}}$ for any $1 \leq \alpha \leq n).\;$
Notice $D$ is just defined locally. Define
a new operator $\cal {L}$ as ( $\phi$ is any test function):
\[   {\cal {L} } \phi = \displaystyle \sum_{\alpha,\beta  = 1}^{n+1}
 \;g'^{\alpha \overline{\beta}}  {{\partial^2 \phi}\over
{\partial z_{\alpha} \partial \overline{z}_{\beta}}}             \]
where $(g'^{\alpha \overline{\beta}}) = (g'_{\alpha \overline{\beta}})^{-1} =
\left( g_{\alpha \overline{\beta}} +  {{\partial^2 \varphi}\over
{\partial z_{\alpha} \partial \overline{z}_{\beta}}}\right)^{-1}.\;$
Differentiating both side of equation (\ref{eq:euler3}) by $D,$ we get
\[
      {\cal {L} }\; D \varphi = D \ln f + \displaystyle \sum_{\alpha,\beta = 1}^{n+1} g'^{\alpha \overline{\beta}} D g_{\alpha \overline{\beta}}.
\]
Thus there exists a constant $C$ which depends only on $(V \times {\bf R},g)$
such that
\begin{equation}
  {\cal {L} } D ( \varphi - \varphi_0) \leq C ( 1 + \displaystyle \sum_{\alpha=1}^{n+1} g'^{\alpha \overline {\alpha}})
\end{equation}

  We will now employ a barrier function of the form
\begin{equation}
\nu = (\varphi - \varphi_0) + s\; ( h - \varphi_0) - N \cdot x^2
\end{equation}
near the boundary point, and $s, N$ are positive constants
to be determined. We may take $\delta $ small enough
so that $x$ is small in $\Omega_{\delta} = (V\times {\bf R}) \cap B_{\delta} (0).\;$ The main essence of the proof is:
\begin{lem} For $N$ sufficiently large and $s, \delta$ sufficiently
small, we have
\[
{\cal {L}}\; \nu \leq -{\epsilon \over 4} ( 1 +  \displaystyle \sum_{\alpha=1}^{n+1}
g'^{\alpha \overline {\alpha}}) \;\;{\rm in}\;\; \Omega_{\delta},\;\;
\nu \geq 0\;\; {\rm on} \; \partial \Omega_{\delta}.
\]
\end{lem}
{\bf Proof} Since $ g_{\alpha \overline{\beta}} +  {{\partial^2 \varphi_0}\over
{\partial z_{\alpha} \partial \overline{z}_{\beta}}} \geq \epsilon \delta_{\alpha \overline{\beta}},$ we have
\[
{\cal {L}} (\varphi - \varphi_0) = \displaystyle \sum_{\alpha,\beta=1}^{n+1} 
g'^{\alpha \overline{\beta}} [ (g_{\alpha \overline{\beta}} +  {{\partial^2 \varphi}\over
{\partial z_{\alpha} \partial \overline{z}_{\beta}}}) -( g_{\alpha \overline{\beta}} +  {{\partial^2 \varphi_0}\over
{\partial z_{\alpha} \partial \overline{z}_{\beta}}}) ]  \leq n+1 - \epsilon \;\displaystyle \sum_{\alpha=1}^{n+1}
g'^{\alpha \overline {\alpha}}
\]
and
\[
{\cal {L}} (h - \varphi_0) \leq C_1 ( 1 + \displaystyle \sum_{\alpha=1}^{n+1}
g'^{\alpha \overline {\alpha}})
\]
for some constant $C_1.\;$ Furthermore, ${\cal {L}}\; x^2 = 2 g'^{(n+1)\overline{n+1}}.\;$ Thus
\[
\begin{array}{lcl} {\cal {L}} \nu & = & {\cal {L}} (\varphi-\varphi_0) + 
s \cdot {\cal {L}} (h - \varphi_0) - 2\cdot  N\cdot g'^{(n+1)\overline{n+1}}\\
        & \leq & n+1 - \epsilon \displaystyle \sum_{\alpha=1}^{n+1}
g'^{\alpha \overline {\alpha}} + s C_1 + s C_1 \displaystyle \sum_{\alpha=1}^{n+1}
g'^{\alpha \overline {\alpha}} - 2 N  g'^{(n+1)\overline{n+1}}.
\end{array}
\]
Suppose $ 0 < \lambda_1 \leq \lambda_2 \leq \cdots \leq \lambda_{n+1}$
are eigenvalues of $ (g'_{\alpha \overline{\beta}})_{(n+1)(n+1)}.\;$ Thus
\[
\displaystyle \sum_{\alpha=1}^{n+1} g'^{\alpha \overline {\alpha}}
=\displaystyle \sum_{\alpha=1}^{n+1} {\lambda_{\alpha}}^{-1},\qquad g'^{(n+1)\overline{n+1}} \geq \lambda_n^{-1}. \]
Thus,
\[
\begin{array}{lcl} {\epsilon \over 4} \displaystyle \sum_{\alpha=1}^{n+1} g'^{\alpha \overline {\alpha}} + N g'^{(n+1)\overline{n+1}} &\geq & {\epsilon \over 4} \displaystyle \sum_{\alpha=1}^{n} {\lambda_{\alpha}}^{-1} + (N +{\epsilon \over 4}) {\lambda_{n+1}}^{-1} \\
 & \geq & (n+1) {\epsilon \over 4} N^{{1\over {n+1}}} (\lambda_1 \cdot \lambda_2 \cdots \lambda_{n+1})^{-{1 \over {n+1}}} = C_2 N^{ {1 \over {n+1}}}. 
\end{array}
\]
Choose $N$ large enough so that
\[
  -C_2 N^{ {1 \over {n+1}}} + (n+1) + s C_1 < - {\epsilon \over 4}.
\]
Choose $s$ small enough so that $ s \cdot C_1 \leq  {\epsilon \over 4}.\;$
Then
\[
  {\cal {L}} \nu \leq - {\epsilon \over 4} ( 1 + \displaystyle \sum_{\alpha=1}^{n+1} g'^{\alpha \overline {\alpha}})). 
\]
From now on we fix $N.\;$
Observe that $\triangle (h -\varphi_0) < - 2\epsilon,$ then
there exists a constant $C_0$ which depends only on $g$ such
that $ h - \varphi_0 > C_0\; x$ near $\partial (V \times {\bf R}).\;$
Choose $\delta $ small enough so that
\[
  s (h -\varphi_0) - N x^2 \geq (s C_0 - N \delta)x \geq  0.
\]
Then $\nu \geq 0$ in $\partial \Omega_\delta.\;$ QED.\\

\noindent {\bf Proof of Lemma 3}:  Let $M = \displaystyle \max (|\nabla \varphi|_g + 1).$
Choose $ A \gg B \gg C,C_1.\;$
In additional, choose $A, B$ as a big multiple of $M.\;$ Notice
that $ |D \varphi | \leq 2 M $ in $\Omega_{\delta}.\;$ For $\delta$
fixed as in Lemma 4, we have $B \delta^2 - | D(\varphi - \varphi_0)| > 0.\;$
Consider $w = A\; \nu + B \;|z|^2 + D(\varphi-\varphi_0).\;$ Then
$w \geq 0$ in $\partial \Omega_{\delta}$ and $ w(0) = 0.\;$ Moreover,
\[
  {\cal {L}} w \leq ( -{{\epsilon A }\over {4}} + 2 B + C) ( 1 + \displaystyle \sum_{\alpha=1}^{n+1} g'^{\alpha \overline {\alpha}}) < 0.
\]  
Maximal Principal implies that $w \geq 0$ in $\Omega_{\delta}.\;$
Since $w(0) = 0,$ then $ {{\partial w}\over {\partial x}} \geq 0.\;$
In other words,
\[
  {{\partial} \over {\partial x}} D \varphi (0) <  C_3 \cdot M
\] 
for some uniform constant $C_3.\; $ Since $D$ is any
1st order constant operator near $\partial (V \times {\bf R}).\;$
Replace $D$ with $-D,$ we get
 \[
  - {{\partial} \over {\partial x}} D \varphi (0) <  C_3 \cdot M
\] 
On the other hand, since $\partial {\bf R}$ is given by $x=0$
in our special case, we then have the trivial estimate:
\[
  {{\partial} \over {\partial y}} D (\varphi-\varphi_0) (0) = 0. 
\]
Therefore,
\[
  | {{\partial} \over {\partial z_{n+1}}} D \varphi (0) | <  C_3 \cdot M
\]
Lemma 3 follows from here directly. QED.\\

\subsection {Blowing up analysis}
\begin{lem} Any bounded weakly sub-harmonic function in two
dimensional plane is a constant.
\end{lem}
This is a standard fact in geometry analysis, we will omit
the proof here. Notice this lemma is false if dimension is no
less than 3.

The essence of blowing up analysis is to use ``micro-scope''
to analyze what happen in a small neighborhood via rescaling.
Hence it doesn't make any difference what the global
structure of background metric is, or what the  metric is. Under rescaling,
everything become Euclidean anyway. We may as well view the manifold
as a domain in Euclidean space. we will use variable $x$ to
denote position in $V \times {\bf R}.\;$
 
\noindent {\bf Proof of theorem 2}: Suppose ${1\over \epsilon_i} = 
\displaystyle \max_{V \times {\bf R}} |\nabla \varphi_i|_g \rightarrow \infty.\;$ We want to draw a contradiction from this statement.

Suppose $ |\nabla \varphi_i|_g(x_i) = {1\over \epsilon_i}.\;$ By theorem
1, we have $\displaystyle \max_{V \times {\bf R}}  \triangle \varphi_i \leq 
{1\over \epsilon_i^2}.\;$ Choose a convergent subsequence of $x_i$
such that $x_i \rightarrow \underline{x}.\;$  Choose a tiny neighborhood $B_{\delta}(\underline{x})$
of $\underline{x} $ so that $g_{\alpha\overline{\beta}}(\underline{x}) = \delta_{\alpha\overline{\beta}}$ and $g $ is essentially an identical matrix in  $B_{\delta}(\underline{x}).\;$For simplicity, let us pretend that $g$ is an Euclidean metric in $B_{\delta}(\underline{x}).\;$
There are two cases to consider:
the first case is when $\underline{x} \in \partial (V\times {\bf R})$ and the
2nd case is when $\underline{x}$ is in the interior of $V\times {\bf R}.\;$\\

We define the blowing up sequence as
\[
   \tilde{\varphi}_i(x) = \varphi_i (x_i +\epsilon_i x), \forall x \in B_{{\delta\over \epsilon_i}}(0).
\]

Then $ |\nabla \tilde{\varphi}_i(0) =1 $ and
\[
 \displaystyle \max_{B_{{\delta\over \epsilon_i}}(0)} |\nabla \tilde{\varphi}_i | \leq 1, \qquad {\rm and}\; \displaystyle \max_{B_{{\delta\over \epsilon_i}}(0)} 
|\triangle \tilde{\varphi}_i | \leq C.
\]

Observe  $\varphi_0 \leq \varphi_i \leq h \;(\;\forall i).\;$ Re-scale
$\varphi_0$ and $h$ accordingly:
\[
   \tilde{\varphi}_0 (x) = \varphi_0 (x_i + \epsilon_i x),\qquad
\tilde{h}(x) = h(x_i + \epsilon_i x),\;\; \forall \; x\; \in B_{{\delta\over \epsilon_i}}(0).
\]

Thus $\displaystyle \lim_{i \rightarrow \infty} \tilde{\varphi}_0 (x)= \varphi_0(\underline{x})$ and $ \displaystyle \lim_{i \rightarrow \infty} \tilde{h} (x)= h(\underline{x}).\;$ Moreover,
\begin{equation}
\tilde{\varphi}_0 \leq \tilde{\varphi}_i \leq \tilde{h},\qquad \forall i = 1, 2,\cdots. \label{eq:c0bound}
\end{equation}

There exists a subsequence of $\tilde{\varphi}_i$ and a limit
function $ \tilde{\varphi}$ in $C^{n+1}$ (or half plane
in case $\underline{x}$ in the boundary) such that in any fixed ball of
$B_{l}(0)$(or half ball if $\underline{x}$ is in the boundary)
we have $\tilde{\varphi}_i \rightarrow \tilde{\varphi}$
in $C^{1,\eta}$ in the ball $ B_{l}(0)$ (or half ball) for 
any $0 < \eta < 1.$  This implies 
\begin{equation}
 |\nabla \tilde{\varphi}(0)| = 1.
\label{eq:gradient=1}
\end{equation}

In additional, inequality (\ref{eq:c0bound}) holds in the limit:
\begin{equation}
\varphi_0(\underline{x}) \leq \tilde{\varphi}(x) \leq h(\underline{x}),\qquad \forall\; x.
\label{eq:limitc0bound}
\end{equation}

Case 1: Suppose $\underline{x} \in \partial (V\times {\bf R}).\;$
Then $h(\underline{x}) = \varphi_0(\underline{x}).\;$  Inequality (\ref{eq:limitc0bound} )
implies that $\tilde{\varphi}$ is a constant function in its domain.
In particular, we have $|\nabla \tilde{\varphi}(x)| \equiv 0.\;$
This  contradicts our assertion (\ref{eq:gradient=1}). Thus the theorem is proved in
this case. \\

Case 2: Suppose $\underline{x}$ is in the interior of $ V \times {\bf R}.\;$
Then $\tilde{\varphi}(x)$ is a well defined $C^{1,\eta}$ and bounded function 
 in $C^{n+1}.\;$ We claim that this
function is weakly sub-harmonic in any complex line through
origin. If this claim is true, then Lemma 5 says it must be
constant for any complex line through origin. Therefore, the
function itself must be a constant as well. Thus $|\nabla \tilde{\varphi}|
\equiv 0.\;$ It again contradicts with our assertion (\ref{eq:gradient=1}). Thus the
theorem is proved also, provided  we can prove this 
claim. \\

  Without loss of generality, we consider the complex line $T$
is
\[
   z_2 = z_3 =\cdots =z_{n+1} =0.
\]

 Observed that (near $\underline{x}$) the following holds
\[
    0 < (\delta_{\alpha\overline{\beta}} + {{\partial ^2 \varphi_i}\over{\partial z_{\alpha} \partial \overline{z}_\beta}})_{(n+1) (n+1)} < { C \over {\epsilon_i^2}} (  \delta_{\alpha\overline{\beta}}  )_{(n+1) (n+1)}, \forall \;i. 
\]
After rescaling, we have
\[
 0 < \epsilon_i^2 \cdot (\delta_{\alpha\overline{\beta}})_{(n+1) (n+1)} + ({{\partial ^2 \tilde{\varphi}_i}\over{\partial z_{\alpha} \partial \overline{z}_\beta}})_{(n+1) (n+1)} < C \cdot (  \delta_{\alpha\overline{\beta}}  )_{(n+1) (n+1)}. 
\]
 Restricting this to a complex line $T,$ we have
\begin{equation}
 0 < \epsilon_i^2 + {{\partial ^2 \tilde{\varphi}_i}\over{\partial z_{1} \partial \overline{z}_1}} < C
\end{equation}
Thus one can choose a subsequence of $\tilde{\varphi}_i$ which
converges $C^{1,\eta} (0 < \eta < 1)$ locally in $T$ to some function
$\psi.\;$ Since the convergence is in $C^{1,\eta},$
thus $\psi = \tilde{\varphi}|_{T};$ i.e., $\psi $ is the restriction
of $\tilde{\varphi}$ in this complex line $T.\;$
By taking weak limit in inequality (11), then
$\tilde{\varphi}_i |_{T} $     weakly converge to
 $\psi $
 in $H^{2,p}_{loc}$ topology for any $p > 1.\;$ Therefore, $\psi $
 is a weakly sub-harmonic
function by taking weak limit in inequality (11). Therefore
$\psi = \tilde{\varphi}|_{T}$ is a constant by
Lemma 5. Our claim is then proved.  QED.\\

\section{Uniqueness of weak $C^0 $ geodesic}
Notation follows from previous section.
\begin{defi} A function $\varphi$ is {\it generalized-pluri-subharmonic} 
in $V\times {\bf R}$ if $\displaystyle \sum_{\alpha,\beta = 1}^{n+1} (g_{\alpha \overline{\beta}} + {{\partial^2 \varphi}\over{\partial z_{\alpha}  \partial \overline{z}_{\beta}}}) dz_{\alpha} d\,\overline{z}_{\beta} $
defines a strictly positive K\"{a}hler metric in $V \times {\bf R}.\;$
\end{defi}\

\begin{defi} A continuous function $\varphi$ in $V\times {\bf R}$ is
a  weak $C^0$  solution to degenerated Monge-Ampere equation (\ref{eq:euler1})
with prescribing boundary data $\varphi_0$ if the following
statement is true: $\forall \; \epsilon > 0,$ there exists a pluri subharmonic
function $\tilde{\varphi}$ in $V\times {\bf R}$ such that $ |\varphi -\tilde{ \varphi}| < \epsilon$ and $\tilde{\varphi}$ solves equation (\ref{eq:euler3}) with 
some positive function $0 < f < \epsilon $ at $t=1,\; $ and with
the same boundary data $\varphi_0.\;$
\end{defi}
Clearly, the solution we obtain through continuous method is
a weak $C^0$ solution of equation (\ref{eq:euler1}).

\begin{theo} Suppose $\varphi_1,\varphi_2$ are two $C^0$ weak solutions
to the degenerated Monge-Ampere equation with prescribing boundary condition
$h_1, h_2.\;$ Then 
\[
       \displaystyle \max_{V\times {\bf R}}\; |\varphi_1 -\varphi_2| \leq
 \max_{\partial (V\times {\bf R})} |h_1 -h_2|.
\]
\end{theo}

\begin{cor} The solution to degenerated Monge-Ampere equation is unique
as soon as the boundary data is fixed.
\end{cor}
{\bf Proof}: Suppose $\phi_1,\phi_2$ are two approximate {\it generalized-pluri-subharmonic}
solutions  of $\varphi_1,\varphi_2$ in the sense of definition 2.
In other words
\[
   det\;(g + {{\partial^2 \phi_i}\over{\partial z_{\alpha}  \partial \overline{z}_{\beta}}}) = f_i\cdot  det\;(g) > 0  \;\; {\rm in}\; V\times {\bf R}; \qquad {\rm and}\;\;
\phi_i = h_i \;{\rm in}\;\; \partial (V\times {\bf R}),\;\; i=1,2 
\]
such that $ \displaystyle \max_{V\times {\bf R}} 
(\;|\varphi_1 - \phi_1| + f_1)$
and $\displaystyle \max_{V\times {\bf R}} (\;|\varphi_2 - \phi_2| + f_2)$
could be made as small as we wanted.\\

$\forall \epsilon > 0,$ we want to show 
\[
   \displaystyle \max_{V\times {\bf R}}\; ( \varphi_1 - \varphi_2) \leq \displaystyle \max_{V\times {\bf R}}\;(h_1 -h_2) + 2 \epsilon.\] 
Choose $f_1$ such that $0 < f_1 < \epsilon$ and
$ \displaystyle \max_{V\times {\bf R}} \;|\varphi_1 - \phi_1|  < \epsilon.\;$
Choose $f_2$ such that $0 < f_2 \leq {1\over 2} \displaystyle \min_{V\times {\bf R}}\; f_1 < \epsilon $ and $ \displaystyle \max_{V\times {\bf R}} \;|\varphi_2 - \phi_2|  < \epsilon.\;$ Then $\phi_1$ is a sub-solution to $\phi_2$
(thus $\phi_1 < \phi_2$) if $h_1 = h_2.\;$ In general, we have
\[
   \displaystyle \max_{V\times {\bf R}}\; (\phi_1 - \phi_2) \leq 
\displaystyle \max_{\partial (V\times {\bf R})}\;(h_1 - h_2).
\]
Thus 
\[
  \begin{array} {lcl} \displaystyle \max_{V\times {\bf R}}\; (\varphi_1 -\varphi_2) & =  &\displaystyle \max_{V\times {\bf R}}\; (\varphi_1 -\phi_1) +\displaystyle \max_{V\times {\bf R}}\; (\phi_1 -\phi_2) + \displaystyle \max_{V\times {\bf R}}\;(\phi_2 - \varphi_2) \\
 & \leq  & \epsilon + \displaystyle \max_{\partial (V\times {\bf R})}\;(h_1 - h_2) +\epsilon \\
  & = & \displaystyle \max_{\partial (V\times {\bf R})}\;(h_1 - h_2) + 2 \epsilon.
\end{array}
\]
Change the role of $\varphi_1$ and $\varphi_2,$ we obtain
\[
\displaystyle \max_{V\times {\bf R}}\; (\varphi_2 -\varphi_1)  \leq \displaystyle \max_{\partial (V\times {\bf R})}\;(h_2 - h_1) + 2 \epsilon.
\]
Thus
\[
\displaystyle \max_{V\times {\bf R}}\; |\varphi_1 -\varphi_2|  \leq \displaystyle \max_{\partial (V\times {\bf R})}\;|h_1 - h_2| + 2 \epsilon.
\]
Let $\epsilon \rightarrow 0, $ we obtain the desired result. QED.

\section{The space of K\"{a}hler metric is a metric space---Triangular inequality}
In this section, we want to prove that the space of K\"{a}hler
metric is a metric space and the $C^{1,1}$ geodesic
between any two points  realizes the global minimal 
length over all possible paths.
To prove this claim, one inevitably  need to take derivatives
of lengths for a family of $C^{1,1}$ geodesics. However,
the length for a $C^{1,1}$ geodesic is just barely  defined 
(the integrand is in $L^p$
space). In general, one can not take derivatives. Therefore, we 
must find ways to circumvent this trouble.

\begin{defi} A path $\varphi(t) (0< t < 1)$ in the 
space of K\"{a}hler metrics is a convex path if $\varphi(t)$
is a {\it generalized-pluri-subharmonic} function in $V\times (I\times S^1)$ (see definition 1).
\end{defi}

Suppose $vol(t) (0\leq t \leq 1)$ is a family of strictly positive volume form
in $V$  such that 
\[
  \displaystyle \int_{V} vol(t) = \displaystyle \int_{V} det\; g.
\]
The notion of $\epsilon$-approximate geodesic is defined
with respect to such a volume form:
\begin{defi}
A convex path $\varphi(t)$ in the space of K\"{a}hler
metrics is called $\epsilon$-approximate 
geodesic if the following holds:
\[
   (\varphi'' - |\nabla \varphi'|_{g(t)}^2)\; det\; g(t) = \epsilon \cdot vol(t)
\]
where $g(t)_{\alpha\overline{\beta}} = g_{\alpha\overline{\beta}} + {{\partial^2 \varphi}\over{\partial z_{\alpha}\partial \overline{z}_{\beta}}}\; (1\leq \alpha,\beta\leq n).\;$

\end{defi}
\begin{rem} The definition is really independent
of these volume forms since we only care what happens when $\epsilon$ is really small. For convenience, sometimes we choose
$vol(t) \equiv det\; g $ 
 (a volume form independent of $t$). 
\end{rem}
\begin{lem} Suppose $\varphi(t) (0\leq t \leq 1)$ is an $\epsilon$-approximate
geodesics. Define the energy element as $E(t) = \displaystyle \int_{V} \varphi'(t)^2 d\;g(t).\;$ Then 
\[
\displaystyle \max_{t}  |{{d \, E}\over {d\,t}}| \leq 2 \;\epsilon \cdot 
\displaystyle \max_{V\times I}\; |\varphi'(t)| \cdot  M \]
where $ M = \displaystyle \int_{V} det \,g$ is the total volume of $V$ which
depends only on the K\"{a}hler class.
\end{lem}
{\bf Proof}: \[
\begin{array} {lcl} |{{d \, E}\over {d\,t}}| & = & |\displaystyle \int_V ( 2\varphi''\varphi' +\varphi'^2 \triangle_{g(t)} \varphi') \; d\, g(t)| \\
   & = & 2 \; |\displaystyle \int_V  \varphi'( \varphi'' - {1\over 2} |\nabla \varphi'|_{g(t)}^2) det \, g(t)|\\ & = & 2\;|\displaystyle \int_V \varphi'\, \epsilon\, vol(t)|
   \leq   2 \; \epsilon \cdot   \displaystyle \max_{V\times I}\; |\varphi'(t)| \cdot  M. \end{array}\qquad {\rm QED}.
\]   

\begin{prop} Suppose $\varphi(t)$ is a $C^{1,1}$ geodesic in $\cal H$ from
 $0$ to $\varphi$ and $I(\varphi)=0.\;$  Then the following inequality holds
\[
\int_{0}^{1} \sqrt{\int_{V} \varphi'^{2}  
d\mu_{{\varphi}_{t}} } dt \geq M^{-1} \left( \max( \int_{\varphi>0} \varphi 
d\mu_{\varphi}, -\int_{\varphi<0} \varphi d\mu_{0}) \right).
\]
In other words, the length of any $C^{1,1}$ geodesic is strictly
positive.
\end{prop}

{\bf Proof}: As in definition (4), suppose $\varphi(\epsilon,t)$
is a $\epsilon-$ approximated geodesic between $0$ and $\varphi.\;$
(We will drop the dependence of $\epsilon$ in this proof since
no confusion shall arise from this omition). First of all, from definition
of $\epsilon-$approximated geodesic, we have

\[  \varphi''- {1\over 2} |\nabla \varphi'|^2_{g(t)} > 0.\]
In particular, we have $\varphi''(t) \geq 0.\;$ Thus
\begin{equation}
   \varphi'(0) \leq \varphi \leq \varphi'(1).
\label{eq:convex1}
\end{equation}
Consider $f(t) = I(t\varphi), t\in [0,1].\;$
Then $f'(t) = \displaystyle \int_V\; \varphi\; d\mu_{t\varphi}$ and
\[
    f''(t) = \displaystyle \int_V \varphi\; \triangle_{g(t\varphi)}\; \varphi\; d\mu_{t\varphi} \leq 0.
\]
Thus, we have $f'(0) \geq { {f(1)-f(0)}\over {1-0}} \geq f'(1).\;$
In other words, we have
\[
  \displaystyle \int_V \varphi\; d\;\mu_0 \geq I(\varphi) \geq \int_V \varphi\; d\; \mu_{\varphi}.
\]

Since we assume $I(\varphi)=0,$ and $\varphi$ not identically zero,
 then it  must take both  positive and negative 
values.  Then the length (or energy) of the geodesic is given by
$$ E=\int_{V}  {\varphi'}^{2} d\mu_{\varphi_{t}}, $$
for any $t\in[0,1]$. In particular, taking $t=1$,
$$ \sqrt{E(1)} \geq  M^{-1/2} \int_{V} \vert \varphi'(1)  \vert d\mu_{\varphi} > M^{-
1/2} \int_{\varphi'(1)>0} \varphi'(1) d\mu_{\varphi}, $$
where $M$ is the volume of $V$ (which is of course the same for all metrics in 
${\cal H}$).                                   It follows
from inequality (\ref{eq:convex1}) that
$$   \int_{ \varphi'(1)>0}  \varphi' \; d\mu_{\varphi}\geq \int_{\varphi>0} \varphi\; 
d\mu_{\varphi}, $$ where the last term is strictly positive by the remarks above, 
and depends only on $\varphi$ and not on the geodesic.  A similar argument gives
    $$   \sqrt {E(0)}> -M^{-1/2} \int_{\varphi<0} \varphi \;d\mu_{0}. $$
The previous lemma implies that for any $t_1,t_2 \in [0,1],$ we have
\[  |E(t_1) - E(t_2)| < C\cdot \epsilon
\]
for some constant $C$ independent of $\epsilon.\;$ Thus
\[
\sqrt{E(t)}  \geq M^{-
1/2} \displaystyle \max (\int_{\varphi>0} \varphi  d\mu_{\varphi},  - \int_{\varphi<0} \varphi d\mu_{0}) - C\cdot \epsilon.
\]
Now integrating from $t=0$ to $1$ and let $\epsilon \rightarrow 0.\;$Then
\[
\int_{0}^1 \sqrt{\displaystyle \int_V \varphi'^2 d\mu_{\varphi}}\;d\;t 
\geq M^{-
1/2} \displaystyle \max (\int_{\varphi>0} \varphi  d\mu_{\varphi}, - \int_{\varphi<0} \varphi d\mu_{0}).
\]
Then this proposition is proved. $QED.$

\begin{rem} This proposition verifies Donaldson's 2nd conjecture \footnote{
In \cite{Dona96}, Donaldson
provided a formal proof to this proposition after assuming the existence of a smooth
geodesic between any two metrics. Our proof follows his idea closely. }. However,
it will not  imply $\cal H$ is a metric space automatically since the
geodesic is not sufficiently differentiable. However, one can easily verifies
that $C^{1,1}$ geodesic minimizes length over all possible convex
curves between the two end points. To show that it minimizes length
 over all possible curves, not just convex ones, we need to 
prove that the triangular inequality is satisfied 
by the geodesic distance (see definition below).
\end{rem} 

\begin{defi} Let $\varphi_1, \varphi_2$ be two distinct
points in the space of metrics. According to theorem 3 and 
Corollary 2, there exists a unique geodesic connecting these
two points. Define the geodesic distance as the length of this geodesic.
Denoted as $d(\varphi_1,\varphi_2).\;$
\end{defi}

\begin{theo}
Suppose $C: \varphi(s): [0,1] \rightarrow {\cal {H}}$ is a smooth
curve in $\cal {H}.\;$ Suppose $p$ is a base point of $\cal {H}.\;$
 For any $s,$ the geodesic distance from
$p$ to $\varphi(s)$ is no greater than the sum of geodesic
distance from $p$ to $\varphi(0)$ and the length from $\varphi(0)$
to $\varphi(s)$ along this curve $C.\;$ In particular, if $C: \varphi(s): [0,1] \rightarrow {\cal {H}}$ is a geodesic, then the geodesic distance
satisfies:
\[
 d(0,\varphi(1)) \leq d(0,\varphi(0)) + d(\varphi(0),\varphi(1)).
\]
\end{theo}

\begin{lem} (Geodesic approximation lemma): Suppose $C_i: \varphi_i(s):
[0,1] \rightarrow {\cal {H}} (i = 1, 2)$ are two smooth
curves in ${\cal {H}}.\;$ For $\epsilon_0$ small enough,
there exist two parameters smooth families
of curves $ C(s,\epsilon): \phi(t,s,\epsilon): [0,1]\times [0,1]\times (0,\epsilon_0] ( 0\leq t, s \leq 1, 0 < \epsilon \leq \epsilon_0) $ such that the following properties hold:
\begin{enumerate}
\item For any fixed $s $ and $\epsilon, C(s,\epsilon)$ is a $\epsilon$-approximate geodesic from $\varphi_1(s)$ to $\varphi_2(s).\;$ More precisely,
$\phi(z,t,s,\epsilon) $ solves the corresponding Monge-Ampere equation:
\begin{equation}
    det\;(g + {{\partial^2 \phi}\over{\partial z_{\alpha}  \partial \overline{z}_{\beta}}}) = \epsilon \cdot  det\;(g) \;\; {\rm in}\; V\times {\bf R}; \qquad {\rm and}\;\;
\phi(z',0,s,\epsilon) = \varphi_1(z',s),\;\phi(z',1,s,\epsilon) = \varphi_2(z',s).
\end{equation} 
Here we follows notation in section 3, and $z_{n+1} = t + \sqrt{-1}\, \theta$ where depends of $\phi$ on $\theta $ is trivial.
\item There exists a uniform constant $C$ (which depends only on $\varphi_1, \varphi_2$ such that
\[
   |\phi| + |{{\partial \phi}\over {\partial s}} | + |{{\partial \phi}\over {\partial t}} | < C; \qquad  0 \leq {{\partial^2 \phi}\over {\partial t^2}}  < C,
\qquad {{\partial^2 \phi}\over {\partial s^2}} < C.\]
\item For fixed $s,$ let $\epsilon \rightarrow 0,$ the convex
curve $ C(s,\epsilon)$ converges to the unique geodesic between
$\varphi_1(s)$ and $\varphi_2(s)$ in weak $C^{1,1}$ topology. 
\item Define energy element  along $C(s,\epsilon)$ by 
\[E(t,s,\epsilon ) = \displaystyle \int_{V} |{{\partial \phi}\over {\partial t}}|^2 d\;g(t,s,\epsilon) \]
where $g(t,s,\epsilon)$ is the corresponding K\"{a}hler metric
define by $\phi(t,s,\epsilon).\;$ Then there exists a uniform 
constant $C$ such that
 \[
\max_{t,s}  |{{\partial \, E}\over {\partial\,t}}| \leq \epsilon \cdot C \cdot  M. \]
In other words, the energy/length element converges
to a constant along each convex curve if $\epsilon \rightarrow 0.\;$
\end{enumerate}
\end{lem}

{\bf Proof }:  Everything follows from theorem 3 and 4 and lemma 6
except the bound on $ |{{\partial \phi}\over {\partial s}} |$ and a upbound
on $ {{\partial^2 \phi}\over {\partial s^2}}$ which
follow from maximal principal directly since
\[
  {\cal {L}} ({{\partial \phi}\over {\partial s}} ) = 0
\] 
and
\[
  {\cal {L}} ({{\partial^2 \phi}\over {\partial s^2}} ) = tr_{g'} \{ Hess {{\partial \phi}\over {\partial s}}\cdot Hess {{\partial \phi}\over {\partial s}}\}
\geq 0. \qquad {\rm QED.}
\] 

{\bf Proof of theorem 5}: Apply geodesic approximation lemma
with special case that $\varphi_1(s)\equiv p.\;$ We follow notations in
the previous lemma. For $\epsilon_0$ small enough,
there exist two parameters smooth families
of curves $ C(s,\epsilon): \phi(t,s,\epsilon): [0,1]\times [0,1]\times (0,\epsilon_0] ( 0\leq t, s \leq 1, 0 < \epsilon \leq \epsilon_0) $ such that
\[
    det\;(g + {{\partial^2 \phi}\over{\partial z_{\alpha}  \partial \overline{z}_{\beta}}}) = \epsilon \cdot  det\;(g), \;\; {\rm in}\; V\times {\bf R}; \qquad {\rm and}\;\;
\phi(z',0,s,\epsilon) = 0,\;\phi(z',1,s,\epsilon) = \varphi(z',s).
\]
  Denote the length of the curve $\phi(t,s,\epsilon)$
from $p$ to $\varphi(s)$ as $L(s,\epsilon),$ denote
the geodesic distance between $p$ and $\varphi(s) $ as $L(s),\;$
and denote the length from $\varphi(0)$ to $\varphi(s)$
along curve $C$ as $l(s).\;$ Clearly, $l(s) = \displaystyle \int_0^s \sqrt{\displaystyle \int_{V} |{{\partial \varphi}\over {\partial \tau}}|^2 d\;g(\tau) }\; d\,\tau$ where $g(\tau)$ is the K\"{a}hler
metric defined by $\varphi(\tau), $ and
\[
  L(s,\epsilon) = \displaystyle \int_0^1 \sqrt{E(t,s,\epsilon)} \,d\,t = \displaystyle \int_0^1 \sqrt{\displaystyle \int_{V} |{{\partial \phi}\over {\partial t}}|^2 d\;g(t,s,\epsilon) }\; d\,t,\;{\rm and}\; \displaystyle \lim_{\epsilon\rightarrow 0} L(s,\epsilon) = L(s).
\]
Define $F(s,\epsilon) = L(s,\epsilon) + l(s)$ and $F(s) = L(s) + l(s).\;$
What we need to prove is : $F(1) \geq  F(0).\;$ This will be done
if we can show that $F'(s) \geq 0,\;\forall \;s\;\in [0,1].\;$
The last statement would be straightforward if the deformation
of geodesics is $C^1.\;$ Since we don't have it, we need
to take derivatives on $F(s,\epsilon)$ for $\epsilon > 0 $
instead. Notice ${{\partial \phi}\over {\partial s}} = 0 $ at $t= 0$
in the following deduction: 
\[
  \begin{array}{lcl} {{d\, L(s,\epsilon)}\over {d\, s}} & = 
& \displaystyle \int_0^1 {1\over 2} \; E(t,s,\epsilon)^{-{1\over 2}} \displaystyle \displaystyle \int_V \left( 2 {{\partial \phi}\over {\partial t}} {{\partial^2 \phi}\over {\partial t \partial s}} + ({{\partial \phi}\over {\partial t}})^2 \triangle_{g(t,s,\epsilon)} {{\partial \phi}\over {\partial s}} \right) d g(t,s,\epsilon)\; d\,t \\
   & = & \displaystyle \int_0^1 \; E(t,s,\epsilon)^{-{1\over 2}} \{{{\partial} \over {\partial
t }}\; \left( \displaystyle \int_V
{{\partial \phi}\over {\partial t}}\;{{\partial \phi}\over {\partial s}}\; d\,
g(t,s,\epsilon)\right) -  \displaystyle \int_V \;{{\partial \phi}\over {\partial s}}\; ({{\partial^2 \phi}\over {\partial t^2}}- {1\over 2} |\nabla {{\partial \phi}\over {\partial t}} |^2 ) d\,g(t,s,\epsilon) \}\; d\,t\\

& = &  \{E(t,s,\epsilon)^{-{1\over 2}}\;\displaystyle \int_V
{{\partial \phi}\over {\partial t}}\;{{\partial \phi}\over {\partial s}}\; d\,
g(t,s,\epsilon))\}|_0^1 \; 
  - \; \displaystyle \int_0^1 \{ E(t,s,\epsilon)^{-{1\over 2}} \displaystyle \int_V \;{{\partial \phi}\over {\partial s}}\; ({{\partial^2 \phi}\over {\partial t^2}}- {1\over 2} |\nabla {{\partial \phi}\over {\partial t}} |^2 ) d\,g(t,s,\epsilon)\} d\,t\\
&  & \qquad + \; \displaystyle \int_0^1 \{ E(t,s,\epsilon)^{-{3\over 2}}\; \displaystyle \int_V
{{\partial \phi}\over {\partial t}}\;{{\partial \phi}\over {\partial s}}\; d\,
g(t,s,\epsilon))\cdot \displaystyle \int_V \;{{\partial \phi}\over {\partial t}}\; ({{\partial^2 \phi}\over {\partial t^2}}- {1\over 2} |\nabla {{\partial \phi}\over {\partial t}} |^2 ) d\,g(t,s,\epsilon)\} \;d\,t\\
& =  &  \displaystyle \int_V
{{\partial \phi(1,s,\epsilon)}\over {\partial t}}\;{{d\, \varphi}\over {d\, s}}\; d\,
g(s) \cdot \{\displaystyle \int_V
|{{\partial \phi(1,s,\epsilon)}\over {\partial t}}|^2 d\,g(s)\}^{-{1\over 2}} 
  -   \; \displaystyle \int_0^1 \{ E(t,s,\epsilon)^{-{1\over 2}} \; \displaystyle \int_V \;{{\partial \phi}\over {\partial s}}\; \epsilon \cdot det\; g \}\; d\,t\\
 & & \qquad \qquad 
+ \displaystyle \int_0^1 \{ E(t,s,\epsilon)^{-{3\over 2}}\; \displaystyle \int_V
{{\partial \phi}\over {\partial t}}\;{{\partial \phi}\over {\partial s}}\; d\,
g(t,s,\epsilon))\cdot \displaystyle \int_V \;{{\partial \phi}\over {\partial t}}\;  \epsilon \cdot det\; g\}\; d\,t.
\end{array}
\]
Observe that By Schwartz inequality, we have
\[
  {{d\, l(s)}\over {d\,s}} = \sqrt{\displaystyle \int_{V} |{{\partial \varphi}\over {\partial s}}|^2 d\;g(s) } \geq -\; \displaystyle \int_V
{{\partial \phi(1,s,\epsilon)}\over {\partial t}}\;{{d\, \varphi}\over {d\, s}}\; d\,
g(s) \cdot \{\displaystyle \int_V
|{{\partial \phi(1,s,\epsilon)}\over {\partial t}}|^2 d\,g(s)\}^{-{1\over 2}}. 
\]
Observe that $F(s,\epsilon) = L(s,\epsilon) + l(s).\;$ Thus
\[
  {{d\, F(s,\epsilon)}\over {d\,s}} \geq 
 - \displaystyle \int_0^1 \{ E(t,s,\epsilon)^{-{1\over 2}} \; \displaystyle \int_V \;{{\partial \phi}\over {\partial s}}\; \epsilon \cdot det\; g \}\; d\,t
  +  \displaystyle \int_0^1 \{ E(t,s,\epsilon)^{-{3\over 2}}\; \displaystyle \int_V
{{\partial \phi}\over {\partial t}}\;{{\partial \phi}\over {\partial s}}\; d\,
g(t,s,\epsilon))\cdot \displaystyle \int_V \;{{\partial \phi}\over {\partial t}}\;  \epsilon \cdot det\; g\}\; d\,t.
\]
Integrating from $0$ to $s \in (0,1],$ we obtain
\[
\begin{array}{lcl}
  F(s,\epsilon) - F(0,\epsilon) & \geq & - \displaystyle \int_0^s \displaystyle \int_0^1 \{ E(t,\tau,\epsilon)^{-{1\over 2}} \; \displaystyle \int_V \;{{\partial \phi}\over {\partial \tau}}\; \epsilon \cdot det\; g \}\; d\,t\; d\,\tau \\
& & \qquad + \displaystyle \int_0^s \displaystyle \int_0^1 \{ E(t,\tau,\epsilon)^{-{3\over 2}}\; \displaystyle \int_V
{{\partial \phi}\over {\partial t}}\;{{\partial \phi}\over {\partial \tau}}\; d\,
g(t,\tau,\epsilon))\cdot \displaystyle \int_V \;{{\partial \phi}\over {\partial t}}\;  \epsilon \cdot det\; g\}\; d\,t \; d\,\tau\\
  &\geq  & - C \epsilon
\end{array}
\]
for some big constant $C$ depends only on $(V\times {\bf R}, g)$ and
the initial curve $C:\varphi(s):[0,1] \rightarrow {\cal {H}}.\;$ Now 
take limit as $\epsilon \rightarrow 0,$ we have $F(s) \geq F(0).\;$
 In other words, the geodesic distance from
$p$ to $\varphi(s)$ is no greater than the sum of geodesic
distance from $p$ to $\varphi(0)$ and the length from $\varphi(0)$
to $\varphi(s)$ along this curve $C.\;$  QED.

\begin{cor} The geodesic distance between any two metrics
realize the absolute minimum of the lengths over all possible
paths.
\end{cor}
{\bf Proof}: For any smooth curve 
$C: \varphi(s): [0,1] \rightarrow {\cal {H}},$
we want to show that the geodesic distance between the two 
end points $\varphi(0)$ and $\varphi(1)$ is no greater than
the length of $C.\;$ However, this follows directly from Theorem
5 by taking $p = \varphi(1)$ and $s=1.\;$ QED.

\begin{theo}  For any two K\"{a}hler potentials $\varphi_1, \varphi_2,$
the minimal length $d(\varphi_1,\varphi_2)$
 over all possible paths which connect
these two K\"{a}hler potentials  is strictly positive,
as long as $\varphi_1 \neq \varphi_2.\;$ In other words, $ ({\cal {H}},d)$ is
a metric space. Moreover, the distance function is at least $C^1.\;$
\end{theo}
{\bf Proof} Immediately from Corollary 3 and proposition 2, we imply
that $({\cal H}, d)$ is a metric space.  Now we want to prove 
the differentiability of distance function.
From the Proof of theorem 5, we have
\[ 
\begin{array}{lcl} {{d\, L(s,\epsilon)}\over {d\, s}}  & = &
 \displaystyle \int_V
{{\partial \phi(1,s,\epsilon)}\over {\partial t}}\;{{d\, \varphi}\over {d\, s}}\; d\,
g(s) \cdot \{\displaystyle \int_V
|{{\partial \phi(1,s,\epsilon)}\over {\partial t}}|^2 d\,g(s)\}^{-{1\over 2}} 
  -  \; \displaystyle \int_0^1 \{ E(t,s,\epsilon)^{-{1\over 2}} \; \displaystyle \int_V \;{{\partial \phi}\over {\partial s}}\; \epsilon \cdot det\; g \}\; d\,t\\
 & & \qquad \qquad + \displaystyle \int_0^1 \{ E(t,s,\epsilon)^{-{3\over 2}}\; \displaystyle \int_V
{{\partial \phi}\over {\partial t}}\;{{\partial \phi}\over {\partial s}}\; d\,
g(t,s,\epsilon))\cdot \displaystyle \int_V \;{{\partial \phi}\over {\partial t}}\;  \epsilon \cdot det\; g\}\; d\,t.
\end{array}
\]

Integrating this from $s_1$ to $s_2$ and divide by $s_2-s_1,$ we have
\[
\begin{array}{l}
  |{{L(s_2,\epsilon)-L(s_1,\epsilon)} \over {s_2 - s_1}} - {1\over{s_2-s_1}} 
\displaystyle \int_{s_1}^{s_2}  \displaystyle \int_V
{{\partial \phi(1,s,\epsilon)}\over {\partial t}}\;{{d\, \varphi}\over {d\, s}}\; d\,
g(s) \cdot \{\displaystyle \int_V
|{{\partial \phi(1,s,\epsilon)}\over {\partial t}}|^2 d\,g(s)\}^{-{1\over 2}}\;d\,s| \\
\leq {1\over {s_2-s_1}}  \displaystyle \int_{s_1}^{s_2} \displaystyle \int_0^1 \{ E(t,s,\epsilon)^{-{1\over 2}}\;\displaystyle \int_V \;|{{\partial \phi}\over {\partial s}}|\; \epsilon \cdot det\; g \}\; d\,t\; d\,s \\
 \\
 \qquad \qquad + {1\over {s_2-s_1}}  \displaystyle \int_{s_1}^{s_2} \displaystyle \int_0^1 \{ E(t,s,\epsilon)^{-{3\over 2}}\; \displaystyle \int_V
|{{\partial \phi}\over {\partial t}}|\;|{{\partial \phi}\over {\partial s}}|\; d\,
g(t,s,\epsilon))\cdot \displaystyle \int_V \;|{{\partial \phi}\over {\partial t}}|\;  \epsilon \cdot det\; g\}\; d\,t\; d\,s \\
\leq C \epsilon.
\end{array}
\]
Let $\epsilon \rightarrow 0,$ and then let $s_2 \rightarrow s_1$  we have
\[\begin{array}{lcl} 
   \displaystyle \lim_{s_2 \rightarrow s_1} {{L(s_2)-L(s_1)} \over {s_2 - s_1}} & = & \displaystyle \lim_{s_2 \rightarrow s_1} {1\over{s_2-s_1}} 
\displaystyle \int_{s_1}^{s_2}  \displaystyle \int_V
{{\partial \phi(1,s)}\over {\partial t}}\;{{d\, \varphi}\over {d\, s}}\; d\,
g(s) \cdot \{\displaystyle \int_V
|{{\partial \phi(1,s)}\over {\partial t}}|^2 d\,g(s)\}^{-{1\over 2}}\;d\,s \\
 & = & \displaystyle \int_V
{{\partial \phi(1,s)}\over {\partial t}}\;{{d\, \varphi}\over {d\, s}}\; d\,
g(s) \cdot \{\displaystyle \int_V
|{{\partial \phi(1,s)}\over {\partial t}}|^2 d\,g(s)\}^{-{1\over 2}}.
\end{array}
\]
The distance function $L$ is then a differentiable function. QED.
\section{Application: Uniqueness of Extremal K\"{a}hler metrics if $C_1 (V) < 0$ and $C_1(V) = 0 $}
 In this section, we want to show that if $C_1 (V) < 0,$ or
if   $C_1 (V)= 0,$  then extremal K\"{a}hler metric is unique
in any K\"{a}hler class.  Furthermore, if $C_1(V) \leq 0,$ extremal
K\"{a}hler metric (if existed) realizes the global minimum of Mabuchi energy
functional in any K\"{a}hler class, thus gave a affirmative answer to
a question raised by Tian Gang in this special case. 
\subsection{Uniqueness of c.s.c metric when $C_1(V)=0$ and the lower
bound of Mabuchi energy for $C_1(V) \leq 0$}
We  should now introduce an important operator--- Lichernowicz
operator $\cal D.\;$ For any function $h,$ 
${\cal D} h = h_{,\alpha \beta} dz^{\alpha} \otimes dz^{\beta}.\;$
If ${\cal D} h = 0,$ then $ \uparrow \overline{\partial} h = g^{\alpha \overline{\beta}} {{\partial h}\over {\partial \overline{\beta}}} {{\partial }\over {\partial  {z_{\alpha}}}} 
$ is a holomorphic
vector field.  Now let us  introduce Mabuchi functional. Like
$I_{\rho}, I,$ it is again defined by its derivatives and one should check
it is well defined by verifying the second derivatives is symmetric (we
will leave this to the reader).
 Let $R$ be the scalar curvature
of metric $g= g_0 + \sqrt{-1} \partial \overline{\partial} \varphi $ and $\underline{R}$ be the average scalar curvature in the cohomology class, let $\psi \in T_{\varphi} {\cal H}.\;$ Then the variation of  Mabuchi energy of $g$ at direction
$\psi$ is:
\[
   \delta_{\psi} E = - \int_{V} \; (R - \underline{R})\cdot \psi \; det\; g.
\] 
Along any smooth geodesic $\varphi(t) \in {\cal H}, $ S. Donaldson shows
\[
   {{d^2 E}\over{d\,t^2}} =  \int_V |{\cal D} \varphi'(t)|^2_g \;det\; g.
\]
Using this, Donaldson shows that constant curvature metric is unique
 in each K\"{a}hler class if the smooth geodesic conjecture is true. 
Now we want to prove
  the uniqueness of constant curvature
metric in each K\"{a}hler class when $C_1(V) < 0$ or $C_1(V) = 0,$
despite the fact we have not proven the smooth geodesic conjecture
yet.
\begin{theo} If either $C_1(V) < 0$ or $C_1(V) = 0,$
then the constant curvature metric (if existed) in any K\"{a}hler
class must be unique.
\end{theo}  
{\bf Proof}: Notations follow from section 5. Suppose $\varphi(t) $
is a $\epsilon-$ approximate geodesic. Then
\[
   det \;g\; (\varphi'' - {1\over 2} |\nabla \varphi'|_g^2 ) = \epsilon \cdot det\; h
\]
where $h$ is a given metrics in the K\"{a}hler class such that
$ Ric(h) < -c \; h$ if $ C_1(V) < 0$ and $ Ric(h) \equiv 0 $ if $ C_1(V) = 0.\;$ Let $f = \varphi'' - {1\over 2} |\nabla \varphi'|_g^2 \geq 0.\;$ Then
\begin{equation}
   \nabla \ln { {det\, g}\over {det\, h}} = - \nabla \ln\, f
\label{eq:Mabuchi1}
\end{equation}

and 
\begin{equation}
  {{d}\over {d\, t}} (\int_{V} \varphi'(t) det\,g ) = \int_V\; f\cdot det g =
\epsilon \cdot \int_V \;det\, h.
\label{eq:Mabuchi2}
\end{equation}
  Let $E$ denote the Mabuchi energy functional. Then
\[
  {{d E}\over{d\,t}} = - \int_{V} \; (R - \underline{R})\cdot \psi \; det\; g 
\]
 A direct calculation yields
\begin{equation}
   {{d^2 E}\over{d\,t^2}} =  \int_V |{\cal D} \varphi'(t)|^2_g \;det\; g 
- \int_V (\varphi'' - {1\over 2} |\nabla \varphi'|_g^2 ) \cdot R\; det \, g +
 \epsilon \cdot \underline{R} \cdot \int_V \;det\, h.
\label{eq:Mabuchi3}
\end{equation}
where we already use equation (\ref{eq:Mabuchi2}). Now the 2nd term
in the right hand side of above equation is:
\[
\begin{array}{ccl} - \int_V \, R \cdot  f \, det g &  = & \int_V {\Delta_g} \ln det \, g \cdot f\cdot det\, g\\
& = & \int_V  \;{\Delta_g} {{\ln det \, g}\over{\ln\, det h}}  \cdot f\cdot det\, g + \int_V \; {\Delta_g} \ln det \, h \cdot f\cdot det\, g\\
& = & - \int_V \nabla_g {{\ln det \, g}\over{\ln\, det h}} \cdot \nabla \ln f\, det g \;-\; \int_V tr_g (Ric(h)) \; f\, det\, g \\
& = & \int_V |\nabla f|_g^2 {1\over f} det\,g - \; \int_V tr_g (Ric(h)) \; f\, det\, g.
\end{array}
\]

Thus integrate from $t=0$ to $1,$
\begin{equation}
\int_{V\times I} | {\cal D} \varphi'|_g^2\, det\, g\, d\,t + \int_{V\times I} 
{{|\nabla f|^2}\over{f}} \, det\, g\, d\,t  - \int_{V\times I} tr_g (Ric(h)) \;f\, det\, g\, d\,t = {{d E}\over{d\,t}} \mid_{0}^{1} - \epsilon \underline{R} \cdot \int_V\, det\, h d\,t. 
\end{equation}
If $\varphi(0) $ and $\varphi(1)$ are both of constant scalar curvature metrics,
then
${{d E}\over{d\,t}} \mid_{0}^{1} = 0$ and
\begin{equation}
\int_{V\times I} | {\cal D} \varphi'|_g^2\, det\, g \,d\,t + \int_{V\times I} 
{{|\nabla f|_g^2}\over{f}} \, det\, g \,d\,t-  \int_{V\times I} tr_g (Ric(h)) \;f \, det\, g\, d\, t  = - \epsilon  \underline{R} \cdot \int_V\, det\, h\, d\,t. 
\end{equation}
Observe that $f \det \, g = \epsilon \cdot det h .\;$ We then imply from
previous equation
\[
 \int_{V\times I} {{ | {\cal D} \varphi'|_g^2}\over f}\, det\, h + \int_{V\times I} 
{|\nabla \ln f|_g^2} \, det\, h - \int_{V\times I} tr_g (Ric(h)) \, det\, h =  - \underline{R} \int_V\, det\, h.
\]
Clearly, if $C_1(V) = 0,$ then $\underline{R} = 0.\; $ Thus
\[ \int_{V\times I} {{| {\cal D} \varphi'|_g^2}\over f}\, det\, h d\,t + \int_{V\times I} 
{|\nabla \ln f|_g^2} \, det\, h 
d\,t = 0
\]
This easily implies that $ {\cal D} \varphi(t) \equiv 0 $  and $ \uparrow \overline{\partial} \varphi'(t) $ is a holomorphic vector field. Since $C_1 = 0,$ the
only holomorphic vector field is constant vector field. Thus $\varphi'(t)$ is
constant on $V$ direction. In other words, $\varphi'(t) $  is a functional of $t$ only. Hence, there
exist at most one constant scalar  curvature metric in each K\"{a}hler class
when $C_1 = 0.\;$ We postpone the proof of case $C_1 < 0$ to the next subsection.

\begin{theo} If $C_1(V) \leq 0, $ then constant scalar curvature 
metric, if existed, realizes the global minimum of
Mabuchi energy functional in each K\"{a}hler  class. In other words,
if Mabuchi energy doesn't have a lower bound, then 
there exists no constant curvature metric in that cohomology
class.
\end{theo}
{\bf Proof} Suppose $\varphi_0 \in \cal H $ is a metric
of constant curvature, then 
\[
 {{d E}\over{d\,t}}|_{\varphi_0} =  - \int_{V} \; (R - \underline{R})\cdot \psi \; det\; g = 0
\]
For any metric $\varphi(1),$ let $\varphi(t) (0\leq t \leq 1)$ is a
path in $\cal H$ which connects between $\varphi(0)$ and $\varphi(1).\;$
In additional, let us assume this is a $\epsilon-$ approximate
geodesic where $\epsilon > 0 $ may be chosen arbitrary small.
From equation (17), we have
\[
\begin{array} {ccl}
   {{d^2 E}\over{d\,t^2}} & = &  \int_V |{\cal D} \varphi'(t)|^2_g \;det\; g 
- \int_V (\varphi'' - {1\over 2} |\nabla \varphi'|_g^2 ) \cdot R\; det \, g +
 \epsilon \cdot \underline{R} \cdot \int_V \;det\, h\\
  & = &  \int_V |{\cal D} \varphi'(t)|^2_g \;det\; g  + \int_V |\nabla f|_g^2 {1\over f} det\,g - \; \int_V tr_g (Ric(h)) \; f\, det\, g + \epsilon \cdot \underline{R} \cdot \int_V \;det\, h\\
  & > & - C\epsilon.
\end{array}
\]
The last inequality holds since the average of scalar curvature
is  a topological invariant. Thus
\[
  E(t) -E(0) \geq -C\epsilon {{t^2}\over 2}, \qquad \forall t \in [0,1].
\]
In particular, this holds for $t=1$
\[
  E(\varphi(1)) -  E(\varphi(0)) = E(1) - E(0) \geq - {{ C\cdot \epsilon}\over 2}.
\]
Now let $\epsilon \rightarrow 0,$ we have
\[
  E(\varphi(1)) \geq E(\varphi(0)).
\]
Thus the theorem is proved since $\varphi(1)$ is arbitrary chosen.
\subsection{Uniqueness of c.s.c. metric when $C_1 < 0$}

 Now we turn our attentions to the case $C_1 < 0.\;$ By initial assumption,
$Ric(h) < -c h$ for some positive constant $c > 0.\;$ Thus
\begin{equation}
\int_{V\times I} {{| {\cal D} \varphi'|_g^2}\over{f}} \, det\, h + \int_{V\times I} 
{|\nabla \ln f|_g^2} \, det\, h + c \cdot \int_{V\times I} tr_g (h) \, det\, h \leq C (= -\underline{R} \int_V\, det\, h).
\label{eq:patch0}
\end{equation}
We want to show that in the limit as $\epsilon \rightarrow 0,$
we still have ${\cal D}\varphi'(t) = 0.\;$
Let us first get an integral estimate on $f^{q\over {2-q}} (1 < q < 2)$
 with respect to measure
$det\,h\; d\,t:$
\[\begin{array}{ccl}
\displaystyle \int_{V\times I} f^{q\over {2-q}} det\, h \; d\, t & \leq & 
C\cdot  \displaystyle \int_{V\times I} f det\, h \; d\, t \\
&\leq  & C\cdot  \displaystyle \int_{V\times I} \{ f \cdot { {det\, g}\over {det\, h}} \}^{1\over n} \cdot  \{{ {det\, h}\over {det\, g}} \}^{1\over n}
 det\, h \; d\, t \\
&\leq  & {\epsilon}^{1\over n} \displaystyle \int_{V\times I} \{{ {det\, h}\over {det\, g}} \}^{1\over n}
 det\, h \; d\, t \\
& \leq &  C \cdot {\epsilon}^{1\over n} \displaystyle \int_{V\times I} tr_g (h)
 det\, h \; d\, t \rightarrow 0.
\end{array}
\]

Let $X =  \uparrow \overline{\partial} \varphi'(t) = g^{\alpha \overline{\beta}} {{\partial \varphi'}\over{\partial \overline{z_{\beta}}}} {{\partial}\over{\partial z_{\alpha}}}.\; $ Then we want to show that $X$ is uniformly
in  $L^2$ with respect  to measure $h + d\,t^2.\;$
\[
\begin{array} {ccl}\int_{V\times I} |X|_h^2 \;det\, h \;d\,t & = & \displaystyle \int_{V\times I} \displaystyle \sum_{\alpha,\beta} h_{\alpha \overline{\beta}} X^{\alpha} \overline{X^{\beta}} \;det\, h \;d\,t \\
& = & \displaystyle \int_{V\times I}\; \displaystyle \sum_{\alpha,\beta,\gamma,\delta} h_{\alpha \overline{\beta}}
g^{\alpha \overline{\gamma}} {{\partial \varphi'}\over{\partial \overline{z_{\gamma}}}}  \overline{\{g^{\beta \overline{\delta}} {{\partial \varphi'}\over{\partial \overline{z_{\delta}}}}\}}\; det\, h \;d\,t\\
& = & \displaystyle \int_{V\times I}\;\displaystyle \sum_{\alpha,\beta,\gamma,\delta} h_{\alpha \overline{\beta}} g^{\alpha \overline{\gamma}} g^{\delta \overline{\beta}} {{\partial \varphi'}\over{\partial \overline{ z_{\gamma}}}}  {{\partial \varphi'}\over{\partial {z_{\delta}}}} \; det\, h \;d\,t\\
& \leq & \displaystyle \int_{V\times I} tr_g (h) |\nabla \varphi'|^2_g \; det\,h\\
& \leq &  C \cdot \displaystyle \int_{V\times I} tr_g (h)  \; det\,h \;d\,t\leq C.
\end{array}
\]
The second to last inequality holds since $f = \varphi'' - {1\over 2} |\nabla \varphi'|_g^2 \geq 0$ and $\varphi'' < C.\;$ Thus $X \in L^2( V \times I)$ has a
uniform upbound for the $L^2$ norm. 

Consider $|D\varphi'|_g $ as a function in $L^2 (V \times I). \;$
First of all, it has a weak limit in $L^2 (V\times I);$ secondly, its
$L^{q} (1 < q < 2) $ norm tends to $0$ as $\epsilon \rightarrow 0.$
\[
\begin{array}{ccl}
 \displaystyle \int_{V\times I}  {|{\cal D}\varphi'|_g}^q \; det\, h \; d\, t
& =  &  \cdot \displaystyle \int_{V\times I}  { {|{\cal D}\varphi'|_g}^q \over {f^{l}}} \cdot f^l \; det\, h \; d\, t \\
& \leq &  \cdot ( \displaystyle \int_{V\times I}  { {|{\cal D}\varphi'|_g}^{s q} \over {f^{l s}}} \; det\, h \; d\, t)^{1\over s} \cdot  ( \displaystyle \int_{V\times I}  f^{l \tau } \; det\, h \; d\, t)^{1\over \tau} \qquad ({\rm where}\; {1\over s} + {1\over \tau} = 1).
\end{array}
\]
Now $l$ is some number we should choose appropriately:
\[
  l s = 1; q s = 2;  {1\over s} + {1\over \tau } = 1.
\]
Thus for any $q < 2,$ we have
\[ s= {2 \over q}; l = {q\over 2}; \tau = {2 \over {2-q}}. \]
 
Thus the above inequality reduce to
\[
\displaystyle \int_{V\times I}  {|{\cal D}\varphi'|_g}^q  \; det\, h \; d\, t \leq 
C \cdot (\displaystyle \int_{V\times I}  { {|{\cal D}\varphi'|_g}^{2} \over {f^{}}} \; det\, h \; d\, t)^{2 \over q} \cdot  (\displaystyle \int_{V\times I}  f^{q\over{2-q}} \; det\, h \; d\, t)^{(2-q)\over 2 } \rightarrow 0.
\]
For any vector $Y \in T^{1,0} (V \times I)$ (i.e., $Y = \displaystyle \sum_{i=1}^n \; Y^{i} {{\partial }\over {\partial z_i}}\;$ where $z_1, z_2, \cdots z_n$
are all of the coordinate functions  in a local chart in $V$.  We use ${{\partial Y} \over {\partial \overline{z}}}$ to denote the vector valued (0,1) form
  $ \displaystyle \sum_{i,j=1}^n {{\partial Y^i} \over {\partial \overline {z_j}}} \;{{\partial }\over {\partial z_i}}\otimes
d\;\overline {z_j}.)\;$ For a scalar function $\psi$ in $V\times I,$ denote $ {{\partial \psi} \over {\partial \overline{z}}}$
as $ \displaystyle \sum_{j=1}^n {{\partial \psi } \over {\partial \overline {z_j}}}  d\;\overline {z_j}.\;$Now
the norm of ${{\partial Y} \over {\partial \overline{z}}}$ and $ {{\partial \psi} \over {\partial \overline{z}}}$
in terms of the metric $h$ are:
\begin{equation}
  \mid{{\partial Y}\over {\partial \overline{z}}}\mid_h^2 = \displaystyle \sum_{\alpha,\beta,r,\delta=1}^n\; h_{\alpha \overline{r}} h^{\overline{\beta} \delta} {{\partial Y^{\alpha}}\over {\partial \overline{z_{\beta}}}} \overline{ \left({{\partial Y^{r}}\over {\partial \overline{z_{\delta}}}}\right)} 
\end{equation}
and
\begin{equation}
  \mid{{\partial \psi}\over {\partial \overline{z}}}\mid_h^2 = \displaystyle \sum_{\alpha,\beta=1}^n\;h^{\alpha \overline{\beta}} {{\partial \psi } \over {\partial {z_{\alpha}}}} {{\partial \psi } \over {\partial \overline {z_\beta}}}. 
\end{equation}

We claim the following inequality holds (for some uniform constant $C$):
\begin{equation}
  \mid{{\partial X}\over {\partial \overline{z}}}\mid_h \leq \sqrt{\displaystyle \sum_{\alpha,\beta,r,\delta=1}^n\;h_{\alpha \overline{r}} h^{\overline{\beta} \delta} {{\partial X^{\alpha}}\over {\partial \overline{z_{\beta}}}} \overline{ \left({{\partial X^{r}}\over {\partial \overline{z_{\delta}}}}\right)} }\leq  C \;\sqrt{tr_g (h)}\; |D\varphi'|_g.
\label{eq:patch1}
\end{equation}
This could be proven by choosing a preferred coordinate, where
$h_{i \overline{j}} = \delta_{i \overline{j}} (1 \leq i, j \leq n) $ while 
$g_{i \overline{j}} = \lambda_i  \delta_{i \overline{j}} (1 \leq i, j \leq n)\;$ in an arbitrary point $O.\;$ Here $\lambda_i$ are eigenvalues of metric $g$
in terms of metric $h.\;$ These $\lambda_i$'s are uniformly bounded from
above since $g$ is uniformly bounded from above.  We want to verify the above inequality
in this point $O.\;$ 
\[
\begin{array}{lcl}
 \mid{{\partial X}\over {\partial \overline{z}}}\mid_h^2 & = &
\displaystyle \sum_{\alpha,\beta,a,b=1}^n\;
{{\partial X^{\alpha} } \over {\partial z_{\overline{\beta}}}} 
 {{\partial X^{\overline{a}} } \over {\partial z_{b}}} h_{\alpha \overline{a}} h^{\overline{\beta}\, b} \\
& = & \displaystyle \sum_{\alpha,\beta,a,b,c,d=1}^n\;h_{\alpha \overline{a}} h^{\overline{\beta} \,b} g^{\alpha \overline{c}} {\varphi'}_{,\overline{c} \overline{\beta}} {\varphi'}_{,d\,b} 
g^{\overline{a} d} = \displaystyle \sum_{\alpha,\beta=1}^n\;\delta_{\alpha \overline{a}} \delta^{\overline{\beta} \,b} {1\over {\lambda_{\alpha}}} \delta^{\alpha \overline{c}} {\varphi'}_{,\overline{c} \overline{\beta}} {\varphi'}_{,d\,b}
{1\over {\lambda_{a}}} \delta^{\overline{a} d} \\
 & = & \displaystyle  \sum_{\alpha, \beta=1}^n\;{1\over {\lambda_{\alpha}^2}} {\varphi'}_{,\overline{\alpha} \overline{\beta}} {\varphi'}_{,\alpha \,\beta}
\leq \left(\displaystyle \sum_{\alpha, \beta=1}^n\;  {{\lambda_{\beta}} \over {\lambda_{\alpha}}} \right) \;
\displaystyle \sum_{\alpha, \beta=1}^n\;{1\over {\lambda_{\alpha} \lambda_{\beta}}} {\varphi'}_{,\overline{\alpha} \overline{\beta}}  {\varphi'}_{,\alpha \,\beta}
\\& \leq  & C \cdot tr_g(h) \cdot |D\varphi'|_g^2.
\end{array}
\]
Here $C$ is a uniform constant. From inequality (\ref{eq:patch0}) and
the fact $g$ is bounded from above, we have
\[
\int_{V\times I} \mid \nabla \;log {{det \,g} \over {det\,h}}\mid_h^2
=  \int_{V\times I} \mid \nabla \;log f \mid_h^2 \leq C \cdot \int_{V\times I} \mid \nabla log f \mid_g^2 det \,g d\,t  \leq C.
\]
and
\[
  \int_{V\times I} \left({{det \;h}\over {det\;g}}\right)^{1\over n} \; det \;h \leq \int_{V\times I} \; tr_g (h) \;det\;h \leq C.
\]
From now on, all of the norm, inner product and integration are taken w.r.t.
 metric $h + d\,t^2$ unless otherwise specified. 
Now define a new vector field $Y$ as
\[
   Y = X \cdot {{det \;g}\over {det\;h}}.
\]
Then
\[
  |Y|_h = |X|_h {{det \;g}\over {det\;h}} \leq C.
\]
In other words, $Y$ has uniform $L^{\infty}$ bound. This
implies that  $Y\cdot {{\partial \ln {{det \;g}\over { det\; h}}} \over {\partial \overline{z}}}$ has uniform $L^q$ bound for any $1<q < 2.\;$
Moreover, for any $1< q < 2,$ we have
\[
  \begin{array} {lcl} \displaystyle \int_{V\times I} \mid {{\partial Y}\over {\partial \overline{z}}} - Y \cdot {{\partial \ln {{det \;g}\over { det\; h}}} \over {\partial \overline{z}}}\mid_h^q & = & \displaystyle \int_{V\times I} \left( \mid
{{\partial X}\over {\partial \overline{z}}} \mid_h  {{det \;g}\over {det\;h}}\right)^q \\ & \leq & \displaystyle \int_{V\times I} \;(\sqrt{tr_g (h)}  {{det \;g}\over {det\;h}})^q\; |D\varphi'|_g^q
\\ &\leq & \displaystyle \int_{V\times I} C \;|D\varphi'|_g^q \rightarrow 0.
\end{array}
\]
This immediately implies that ${{\partial Y}\over {\partial \overline{z}}}
$ are uniformly bounded in $L^q$  for any $1<q < 2.\;$\\

Now, all of these quantities, $X, Y, {{\partial Y}\over {\partial \overline{z}}},\;$ and $ {{det \;g}\over {det\;h}},\cdots$ are geometrical
quantities which 
depend on $\epsilon.\;$  Since their respective soblev norms are uniformly controlled,
we can take weak limits of these quantities in some appropriate sense.
Denote the corresponding weak limits (when $\epsilon \rightarrow 0$) as
$X, Y, {{det \;g}\over {det\;h}}, \cdots.\;$ Then $X(\epsilon) \rightharpoonup X$
weakly in $L^2(V \times I),\; Y (\epsilon) \rightharpoonup X$ weakly
in $L^{\infty}(V\times I)\;$ and $ {{\det\;g} \over {det\;h}} (\epsilon)
\rightharpoonup {{det\;g} \over {det\;h}}$ weakly in $L^{\infty} (V\times I),\cdots.\;$ \\

Consider $u = \ln {{det h}\over {det\;g}}.\; $ For simplicity,
assume $ u > 0$ (otherwise $u > -c$ for some positive constants). Then
the following two equations  holds in the limit
\[
   {{\partial Y}\over {\partial \overline{z}}} +  Y \cdot {{\partial u} \over {\partial \overline{z}}} = 0,\qquad {\rm and}\qquad Y= X \;e^{-u}
\]
in the sense of $L^q (V \times I)$ for any $1 < q < 2.\;$ Moreover,  we have the following estimates:
\[
  \int_{V \times I} e^{{1\over n} u} \leq C; \qquad \int_{V \times I} |{{\partial u}\over{\partial \overline{z}}}|^2  \leq C;\qquad {\rm and}\; \int_{V \times I} |X|^2 \leq C.
\]

Now define a new sequence of vectors $X_{,k} (k =1,2,\cdots )$ as
$X_{,k} = Y \displaystyle \sum_{i=0}^k  {{u^i}\over {i!}}.\; $  This
is well defined since $u$ is in $L^p (V \times I) $ for any $p>1.\;$
 Then
\[\begin{array} {lcl}
  |X_{,k}|  & = & |Y| \displaystyle \sum_{i=0}^k  {{u^i}\over {i!}}\\
         & \leq  & (|X| e^{-u}) e^{u} \leq |X|.
\end{array}
\]
The equality holds in the last inequality whenever $e^{-u} \neq 0.\;$
Thus
\[
  \displaystyle \int_{V \times I} \; |X_{,k}|^2 \leq \int_{V \times I} |X|^2 \leq C.
\]
By definition, it is clear $\| X_{,k} \|_{L^2 (V \times I)} \leq 
\| X_{,m} \|_{L^2 (V \times I)}$ whenever $ k \leq m.\;$ Thus, there exists
a positive number $A \leq \| X \|_{L^2 (V \times I)} $ such that 
  $\displaystyle \lim_{ k \rightarrow \infty} \|X_{,k}\|_{L^2 (V \times I)} = A.\; $ For $m > k$, we have
\[
\begin{array} {lcl}
\| X_{,m} \|_{L^2 (V \times I)}^2 & = &  \displaystyle  \int_{V \times I}\; |X_{,m}|^2 \\
  & = & \displaystyle \int_{V \times I} \; |Y|^2  (\displaystyle \sum_{i=0}^m  {{u^i}\over {i!}})^2 \geq
\displaystyle \int_{V \times I} \;  |Y|^2 \left(  (\displaystyle \sum_{i=0}^k  {{u^i}\over {i!}})^2 +  (\displaystyle \sum_{i=k+1}^m  {{u^i}\over {i!}})^2 \right) \\
  & = & \displaystyle \int_{V\times I}\; (|X_{,k}|^2 + |X_{,m} - X_{,k}|^2) = \| X_{,k} \|_{L^2 (V \times I)}^2 + \| X_{,m} - X_{,k} \|_{L^2 (V \times I)}^2. 
\end{array}
\]
Taking limits as $m,k \rightarrow \infty,$ we have
$ \| X_{,m} - X_{,k} \|_{L^2 (V \times I)}^2 \rightarrow 0.\;$
Thus $X_{,k}\; (k =1,2, \cdots )$ is a Cauchy sequence in $L^2 (V \times I)\;$ and there exists
a strong limit $X_{,\infty} $ in $L^2 (V \times I).\;$ By definition,
we know that $X_{,\infty} = X$ almost everywhere in $V \times I$ \footnote{It is easy to prove that $X_{,\infty} = X$ in the sense of $L^q (V\times I) $
for any $ (1<q<2).$}. 
We want to show that $X_{,\infty}$ is weakly holomorphic in  $V-$ direction.
A straightforward calculation yields
\[
  \begin{array}{lcl}  {{\partial X_{,k}}\over {\partial \overline{z}}}  & = & {{\partial Y}\over {\partial \overline{z}}} \displaystyle \sum_{i=0}^k  {{u^i}\over {i!}}
 + Y {{\partial }\over {\partial \overline{z}}} \left(\displaystyle \sum_{i=0}^k  {{u^i}\over {i!}} \right) \\
  & = & - Y {{\partial u}\over {\partial \overline{z}}} \displaystyle \sum_{i=0}^k  {{u^i}\over {i!}} + Y \left(\displaystyle \sum_{i=0}^{k-1}  {{u^i}\over {i!}}\right) \;
{{\partial u}\over {\partial \overline{z}}} \\
  & = & - ( X_{,k} - X_{,k-1}) {{\partial u}\over {\partial \overline{z}}}.
\end{array}
\]
We want to show that ${{\partial X_{,\infty}}\over {\partial \overline{z}}} =0$
in the sense of distribution. We just need to show it in any open set $U \times I$ where $U$ is a coordinate chart in $V.\;$ Denote $(z_1, z_2, \cdots z_n)$
as coordinate variable in $U.\;$ 
Then, for any vector valued smooth function $\psi= (\psi^1,\psi^2,\cdots \psi^n)$
which vanish in $\partial (U\times I),\;$ and for any $1\leq j\leq n.\;$
we have
\[\begin{array} {lcl}
 \mid \displaystyle \int_{V \times I} \; X_{,k} \cdot {{\partial \overline{\psi}}\over {\partial \overline{z_j}}} \mid &  = &  \mid - \displaystyle \int_{V \times I} \;{{\partial X_{,k}}\over {\partial \overline{z_j}}}\cdot \overline{\psi} \mid \\
& = & \mid  \displaystyle \int_{V \times I} \;( X_{,k} - X_{,k-1}) {{\partial u}\over {\partial \overline{z_j}}} \overline{\psi} \mid \\ & \leq &  C \cdot \| X_{,k} -X_{,k-1}\|_{L^2 (V\times I)} \cdot \sqrt {\displaystyle \int_{V\times I} |\nabla u|^2} \\ & \leq & C \| X_{,k} -X_{,k-1}\|_{L^2 (V\times I)}.
\end{array}
\]
Now, taking limit as $k \rightarrow \infty,$ we have
\[
\displaystyle \int_{V \times I} \; X_{,\infty} \cdot {{\partial \overline{\psi}}\over {\partial \overline{z_j}}} = 0,\qquad {\rm for \; any \;} j =1,2,\cdots n 
\]
and for any smooth vector valued function $\psi= (\psi^1,\psi^2,\cdots \psi^n)$ which vanish in $\partial (U \times I).\;$ Thus, $X_{,\infty}$ is a weak holomorphic vector field in $V$ direction for almost all $t$.

 Now recalls that
\[
\int_{V\times I} |X_{,\infty}|_h^2 \;det\, h \;d\,t < C.
\]
This implies that $X_{,\infty}$ is in $L^2(V \times \{t\})$ for almost
all $t \in [0,1].\;$ Since $X_{,\infty}$ is weakly holomorphic in $V \times \{t\}$
for all  $t,$ thus $X_{,\infty}$ must be holomorphic for those $t$ where
$X_{,\infty}$ is in $L^2(V\times \{t\}).\;$  However, there
is no holomorphic vector field in $V$ since $C_1 < 0.\;$ Thus
$X_{,\infty} \equiv 0$ for all of those $t$ where
$X_{,\infty}$ is in $L^2(V\times \{t\}).\;$ This implies that
$X_{,\infty} = 0$ in $V\times I.\;$ Thus $X=0$ since $X = X_{,\infty}$
in the sense of $L^q(V\times I)$ for any $1 < q < 2.\;$ 
Recall
\[
  {{\partial \varphi'(t)}\over {\partial z_{\alpha}}} =  \displaystyle \sum_{\beta=1}^n\; g_{\alpha \overline{\beta}} {{X}}^{\overline{\beta}} =  \displaystyle \sum_{\beta=1}^n\; g_{\alpha \overline{\beta}} {{X_{,\infty}}}^{\overline{\beta}} =  0.
\] 
In other words,  $\varphi'(t)$ is trivial in $V-$direction
 and it is a function of $t$ only for all $t \in [0,1].\; $
 Thus, $\varphi(0) $  and $\varphi(1)$
differ only by a constant in $V$ direction. Therefore they
represent same metric in each K\"{a}hler class.

\bibliography{test}
\nocite{Bourg97}
\nocite{Bourg85}
\nocite{Semmes93}
\nocite{Semmes95}
\nocite{GuanB98}
\nocite{GuanP971}
\nocite{GuanP972}
\nocite{Roch84}
\nocite{Slod96}
\end{document}